\RequirePackage{ifpdf}
\ifpdf 
\documentclass[pdftex]{sigma}
\else
\documentclass{sigma}
\fi

\newtheorem{Theorem}{Theorem}[section]
\newtheorem{Proposition}[Theorem]{Proposition}
\newtheorem{Corollary}[Theorem]{Corollary}

\newtheorem{Conjecture}[Theorem]{Conjecture}
\newtheorem{Remark}[Theorem]{Remark}
\numberwithin{equation}{section}
\newcommand{\ta}{\theta}
\newcommand{\la}{\lambda}

\begin{document}
\allowdisplaybreaks

\renewcommand{\PaperNumber}{056}

\renewcommand{\thefootnote}{$\star$}

\FirstPageHeading

\ShortArticleName{Macdonald Polynomials and Multivariable Series}

\ArticleName{Macdonald Polynomials and\\ Multivariable Basic
Hypergeometric Series\footnote{This paper is a contribution to the
Vadim Kuznetsov Memorial Issue `Integrable Systems and Related
Topics'. The full collection is available at
\href{http://www.emis.de/journals/SIGMA/kuznetsov.html}{http://www.emis.de/journals/SIGMA/kuznetsov.html}}}

\Author{Michael J.\ SCHLOSSER}

\AuthorNameForHeading{M.J.\ Schlosser}

\Address{Fakult\"at f\"ur Mathematik, Universit\"at Wien,
Nordbergstra{\ss}e 15, A-1090 Vienna, Austria}
\Email{\href{mailto:michael.schlosser@univie.ac.at}{michael.schlosser@univie.ac.at}}
\URLaddress{\url{http://www.mat.univie.ac.at/~schlosse/}}

\ArticleDates{Received November 21, 2006; Published online March 30, 2007}

\Abstract{We study Macdonald polynomials from a basic
hypergeometric series point of view. In particular, we show that
the Pieri formula for Macdonald polynomials and its recently
discovered inverse, a recursion formula for Macdonald polynomials,
both represent multivariable extensions of the terminating
very-well-poised ${}_6\phi_5$ summation formula. We derive several
new related identities including multivariate extensions of
Jackson's very-well-poised ${}_8\phi_7$ summation. Motivated by
our basic hypergeometric analysis, we propose an extension of
Macdonald polynomials to Macdonald symmetric functions indexed by
partitions with complex parts. These appear to possess nice
properties.}


\Classification{33D52; 15A09; 33D67}

\renewcommand{\thefootnote}{\arabic{footnote}}
\setcounter{footnote}{0}


\section{Introduction}

The objective of this paper is to study some aspects of $A_{n-1}$
{\em Macdonald polynomials} (which are a family of symmetric
multivariable orthogonal polynomials associated with the
irreducible reduced root system $A_{n-1}$, introduced by
I.G.~Macdonald~\cite{Ma1} in the 1980's), with a particular
emphasis on their connection to (multivariable) {\em basic
hypergeometric series}.

Macdonald polynomials of type $A$ are indexed by integer
partitions, and form a basis of the algebra of symmetric functions
with rational coef\/f\/icients in two parameters $q$ and $t$. They
generalize many classical bases of this algebra, including
monomial, elementary, Schur, Hall--Littlewood, and Jack symmetric
functions. These particular cases correspond to various
specializations of the indeterminates $q$ and $t$. In terms of
basic hypergeometric series, the Macdonald polyno\-mials
correspond to a multivariable generalization of the {\em
continuous $q$-ultraspherical polynomials}, see~\cite{Ko}.

A principal tool for studying $q$-orthogonal polynomials (see
e.g.\ \cite{Is}) is the theory of basic hypergeometric series
(cf.\ \cite{GR}), rich of identities, having applications in
dif\/ferent areas such as combinatorics, number theory,
statistics, and physics (cf.\ \cite{An}). Hypergeometric and basic
hypergeometric series undoubtedly play a prominent role in special
functions, see \cite{AAR}. Even in one variable, they are still an
object of active research. A notable recent advance includes {\em
elliptic} (or {\em modular}) {\em hypergeometric series} (surveyed
in \cite[Ch.~11]{GR} and \cite{Sp}) which is a one-parameter
generalization of basic hypergeometric series, f\/irst introduced
by Frenkel and Turaev~\cite{FT} in a~study related to statistical
mechanics.

A convenient tool (suggested here) for further developing the
theory of multivariable $q$-orthogonal polynomials is the theory
of {\em multivariable basic hypergeometric series associated with
root systems} (or, equivalently, with {\em Lie algebras}). Basic
hypergeometric series associated with the root system $A_{n-1}$
(or equivalently, associated with the unitary group $U(n)$) have
their origin in the work of the three mathematical physicists
Holman, Biedenharn, and Louck, starting in 1976, see
\cite{Ho,HBL}. Their work was done in the context of the quantum
theory of angular momentum, using methods relying on the
representation theory of the unitary group $U(n)$. Subsequently,
extensive investigations in the theory of multiple basic
hypergeometric series associated to the root system $A_{n-1}$ have
been carried out by R.A.~Gustafson, S.C.~Milne, and later various
other researchers. As result, many of the classical formulae for
basic hypergeometric series from \cite{GR} have already been
generalized to the setting of $A_{n-1}$ series (see
Subsection~\ref{secanseries} for some selected results).

An important result that connects $A_{n-1}$ basic hypergeometric
series with Macdonald polynomials is Kajihara and
Noumi's~\cite{KN} explicit construction of raising operators of
row type for Macdonald polynomials. Their construction utilized
$A_{n-1}$ terminating ${}_1\phi_0$ and ${}_2\phi_1$ summations
previously obtained by Milne~\cite{Mi5} (which, however, were
derived independently in \cite{KN} using Macdonald's
$q$-dif\/ference operator).

In this paper we reveal yet more connections of $A_{n-1}$ basic
hypergeometric series with Macdonald polynomials. On the other
hand, we also understand the present work as a contribution
towards the development of a theory of multivariable {\em
very-well-poised}\/ basic hypergeometric series involving
Macdonald polynomials of type $A$.

Various identities for multiple basic hypergeometric series of
Macdonald polynomial {\em argument} have been derived by
Macdonald~\cite[p.~374, Eq.~(4)]{Ma}, Kaneko~\cite{Kan1,Kan2},
Baker and Forrester~\cite{BF}, and Warnaar~\cite{W}. These authors
in fact derived multivariable analogues of many of the classical
summation and transformation formulae for basic hypergeometric
series. As a matter of fact, none of these multivariate identities
reduce to summations or transformations for {\em
very-well-poised}\/ basic hypergeometric series in the univariate
case. There are thus several other classical basic hypergeometric
identities for which higher-dimensional extensions involving
Macdonald polynomials of type $A$ have not yet been explicitly
determined. In this paper we partly remedy this picture by
explicitly pointing out several ``Macdonald polynomial analogues''
of very-well-poised identities. Although some of these identities
(such as the Pieri formula) are not new, their ``very-well-poised
context'' appears so far to have kept unnoticed (at least, not
explicitly mentioned in literature).

The present work makes some further ties between Macdonald
polynomials and multiple basic hypergeometric series associated to
$A_{n-1}$ apparent. In our investigations, we utilize some recent
results from \cite{LS} and application of multidimensional inverse
relations as main ingredients. Since the present paper can be
regarded as a sequel to \cite{LS}, it is appropriate to give a
brief survey here describing those results of \cite{LS} which are
relevant for this paper. This concerns, in particular, the
derivation of a rather general explicit multidimensional matrix
inverse \cite[Th.~2.6]{LS}. The combination of a special case of
this matrix inverse with the Pieri formula for Macdonald
polynomials yielded a new recursion formula for Macdonald
polynomials~\cite[Th.~4.1]{LS}. This formula, the main result of
\cite{LS} -- subsequently also derived by Lassalle \cite{La2} by
alternative means, namely by functional equations -- represents a
recursion for the Macdonald polyno\-mials~$Q_\la$ on the row
length of the indexing partition $\la$. More precisely, it expands
a Macdonald polynomial of row length $n+1$ into products of one
row and $n$ row Macdonald polynomials with explicitly determined
coef\/f\/icients (that are rational functions in $q$ and $t$). The
formula is an $(n+1)$-variable extension of Jing and
J\'{o}zef\/iak's~\cite{JJ} well-known two-variable result and, as
such, the Macdonald polynomial extension of the celebrated
Jacobi--Trudi identity for Schur functions~\cite[p.~41,
Eq.~(3.4)]{Ma}.

In this paper we use special instances of the multidimensional
matrix inverse of \cite[Th.~2.6]{LS} and specif\/ic summation
theorems for multivariable basic hypergeometric series to deduce
several new identities for multivariable basic hypergeometric
series. These turn out to play a special role (or, in some cases,
we believe may be useful) when considering identities involving
principally specialized Macdonald polynomials. In particular, we
establish various new multivariable extensions of the terminating
very-well-poised ${}_8\phi_7$ summation theorem, and of
(terminating and nonterminating) very-well-poised ${}_6\phi_5$
summations. In one case, we provide (only) a conjectural
multivariable terminating very-well-poised ${}_8\phi_7$ summation,
see Conjecture~\ref{cn87n}, which nevertheless we strongly believe
to be true.

We further identify the Pieri formula for Macdonald polynomials
and its inverse, the recursion formula for Macdonald polynomials,
as two dif\/ferent multivariable generalizations of the
terminating very-well-poised ${}_6\phi_5$ summation. More
precisely, this interpretation comes from the explicit structure
the two formulae have after application of analytically continued
principal specialization (coined ``hypergeometric specialization''
in the present context) to both sides of the respective
identities.

Motivated by the basic hypergeometric analysis, we use the
recursion formula to def\/ine~$A_{n-1}$ Macdonald polynomials for
``complex partitions'', no longer indexed by sequences of
non-increa\-sing nonnegative integers but by arbitrary f\/inite
sequences of complex numbers. These $A_{n-1}$ Macdonald symmetric
functions are identif\/ied as certain multivariable nonterminating
basic hypergeometric series which appear to satisfy some nice
properties. We only show a few of them here, in this point our
exposition is mainly intended to hint towards specif\/ic new
objects that we feel deserve future (more thorough and rigorous)
study.

Our paper is organized as follows: In Section~\ref{secmacd} we
recollect some facts on symmetric functions, and, in particular,
$A_{n-1}$ Macdonald polynomials. We conclude this section by
stating the Pieri and the recursion formula, both playing crucial
roles in the paper. In Section~\ref{secbhs} we review some
material on basic hypergeometric series in one and several
variables, hereby explicitly listing several of the summations we
need. In Section~\ref{secmsnt} we derive several identities of a
new type that can be characterized by containing specif\/ic
determinants in the summand. Section~\ref{sechypsp} is devoted to
the hypergeometric specialization of identities involving
$A_{n-1}$ Macdonald polynomials. Section~\ref{secmorehyp} contains
a brief discussion on more basic hypergeometric identities
involving Macdonald polynomials. In Section~\ref{secmacdrev} we
consider Macdonald symmetric functions indexed by partitions with
complex parts. Finally, the Appendix contains the multidimensional
matrix inverse of \cite{LS}, together with a number of useful
corollaries.

\section[Preliminaries on $A_{n-1}$ Macdonald polynomials]{Preliminaries
on $\boldsymbol{A_{n-1}}$ Macdonald polynomials}\label{secmacd}

Standard references for Macdonald polynomials are \cite{Ma1},
\cite[Ch.~6]{Ma}, and \cite{Ma2}. In the following, we recollect
some facts we need.

Let $X=\{x_1,x_2,x_3,\ldots\}$ be an inf\/inite set of
indeterminates and $\mathcal S$ the algebra of symmetric functions
of $X$ with coef\/f\/icients in $\mathbb Q$. There are several
standard algebraic bases of $\mathcal S$. Among these there are
the {\em power sum} symmetric functions, def\/ined by
$p_{k}(X)=\sum\limits_{i \ge 1} x_i^k$, for integer $k\ge 1$. Two
other standard algebraic bases are the {\em elementary} and {\em
complete} symmetric functions $e_{k}(X)$ and $h_{k}(X)$, which are
def\/ined, for integer $k\ge 0$, by their generating functions
\begin{gather*}
\prod_{i \ge 1} (1 +ux_i) =\sum_{k\geq0} u^k\, e_k(X), \qquad
\prod_{i \ge 1}  \frac{1}{1-ux_i} = \sum_{k\geq0} u^k\, h_k(X).
\end{gather*}

A partition $\la= (\la_1,\la_2,\dots)$ is a weakly decreasing
(f\/inite or inf\/inite) sequence of nonnegative integers, with a
f\/inite number of positive integers, called parts. The number of
positive parts is called the length of $\la$ and denoted $l(\la)$.
If $l(\la)=n$, we often suppress any zeros appearing in the
(sequential) representation of $\la$ and write
$\la=(\la_1,\dots,\la_n)$. For any integer $i\geq 1$, $m_i(\la) =
\textrm{card} \{j: \la_j  = i\}$ is the {\em multiplicity} of $i$
in $\la$.  Clearly $l(\la)=\sum\limits_{i\ge1} m_i(\la)$. We shall
also write $\la= (1^{m_1},2^{m_2},3^{m_3},\ldots)$ (where the
parts now appear in increasing order). We set
\begin{gather*}
z_\la  = \prod_{i \ge  1} i^{m_i(\lambda)} m_i(\lambda)!.
\end{gather*}

For any partition $\la$, the symmetric functions $e_{\la}$,
$h_{\la}$ and $p_{\la}$ are def\/ined by
\begin{gather}\label{bas}
f_{\la}=
\prod_{i=1}^{l(\la)}f_{\la_{i}}=\prod_{i\geq1}(f_i)^{m_{i}(\la)},
\end{gather}
where $f_i$ stands for $e_i$, $h_i$, or $p_i$, respectively. (Here
and in the following we sometimes omit wri\-ting out the argument
$X$ of the function, for brevity, assuming there is no confusion.)
These~$e_{\la}$,~$h_{\la}$, and $p_{\la}$ now each form a linear
(vector space) basis of $\mathcal{S}$. Another classical  basis is
formed by the {\em monomial} symmetric functions $m_{\la}$,
def\/ined as the sum of all distinct monomials whose exponent is a
permutation of $\la$.

For any two indeterminates $a$, $q$, we def\/ine the $q$-shifted
factorial as
\begin{gather}\label{ipr}
{(a;q)}_\infty=\prod_{j\ge 0}(1-aq^j),\qquad\text{and}\qquad
(a;q)_k=\frac{(a;q)_{\infty}}{(aq^k;q)_{\infty}},
\end{gather}
for arbitrary $k$ (usually an integer). Sometimes we rather
consider expressions involving the $q$-shifted factorials in
\eqref{ipr} in an analytic (holomorphic) sense. In such a case,
$a$, $q$ are complex numbers, usually with $|q|<1$ (ensuring that
the inf\/inite product in \eqref{ipr} converges, a restriction not
needed when $k$ is an integer and $(a;q)_k$ is indeed a f\/inite
product) and $k$ can be any complex number. We also shall make
frequent use of the condensed notation
\begin{gather}\label{cdn}
(a_1,\dots,a_m;q)_k=(a_1;q)_k\cdots(a_m;q)_k,
\end{gather}
for arbitrary $k$ and nonnegative integer $m$.

Let $\mathbb{Q}(q,t)$ be the algebra of rational functions in the
two  indeterminates $q$, $t$, and
$\mathsf{Sym}=\mathcal{S}\otimes\mathbb{Q}(q,t)$ the algebra of
symmetric functions with coef\/f\/icients in $\mathbb{Q}(q,t)$.

For any $k\ge 0$, the {\em modified complete} symmetric function
$g_{k}(X;q,t)$ is def\/ined by the gene\-rating function
\begin{gather}\label{genserg}
\prod_{i \ge 1} \frac{{(tux_i;q)}_{\infty}}{{(ux_i;q)}_{\infty}}
=\sum_{k\ge0} u^k g_{k}(X;q,t).
\end{gather}
It is often written in plethystic notation \cite[p.~223]{La1},
that is
\begin{gather*}
g_{k}(X;q,t)=h_{k}\left[\frac{1-t}{1-q}\,X\right].
\end{gather*}

The symmetric functions $g_{k}(q,t)$ form an algebraic basis of
$\mathsf{Sym}$. Their explicit expansion in terms of power sums
and monomial symmetric functions has been given by
Macdonald~\cite[pp.~311 and 314]{Ma} and in terms of other
classical bases by Lassalle~\cite[Sec.~10, p.~237]{La1}. The
functions $g_\la(q,t)$, def\/ined as in \eqref{bas} and
\eqref{genserg}, form a linear basis of $\mathsf{Sym}$.

We are ready to def\/ine the Macdonald polynomials
$P_{\la}(X;q,t)$. On one hand, they are of the form (recalling
that $m_\la$ denotes a monomial symmetric function)
\begin{gather*}
\!P_\la(q,t)=m_\la\!+ \text{a linear combination of the $m_\mu\!$
for $\mu$ preceding $\la$ in lexicographical order.}
\end{gather*}
Furthermore, they form an orthogonal basis of $\mathsf{Sym}$ with
respect to the scalar product $\langle\, , \,\rangle_{q,t}$
def\/ined by
\begin{gather*}
\langle p_\la,p_\mu\rangle_{q,t}{}=\delta_{\la \mu}\,z_{\la}\,
\prod_{i=1}^{l(\la)} \frac{1-q^{\la_i}}{1-t^{\la_i}}.
\end{gather*}
Although these two conditions already {\em overdetermine} the
symmetric functions $P_{\la}(X;q,t)$, the latter can be shown to
exist (and are moreover unique), see \cite[p.~322]{Ma}.

Let $Q_{\la}(q,t)$ denote the dual basis of $P_{\la}(q,t)$ for
this scalar product. These are also called Macdonald polynomials,
they dif\/fer from the latter only by a rational function of $q$
and $t$. More precisely, one has
\begin{gather}\label{mcdnf}
Q_{\la}(X;q,t)= b_{\la}(q,t) \,P_{\la}(X;q,t),
\end{gather}
with $b_{\la}(q,t)={\langle
P_{\la}(q,t),P_{\la}(q,t)\rangle_{q,t}}^{-1}$ specif\/ied as
follows (see \cite[p.~339, Eq.~(6.19)]{Ma} and
\cite[Prop.~3.2]{Kan1}):
\begin{gather}
b_{\la}(q,t) =\prod_{1\le i\le j\le l(\la)}
\frac{(q^{\la_i-\la_j}t^{j-i+1};q)_{\la_j-\la_{j+1}}}
{(q^{\la_i-\la_j+1}t^{j-i};q)_{\la_j-\la_{j+1}}} =\prod_{1\le i\le
j\le l(\la)}
\frac{(qt^{j-i};q)_{\la_i-\la_j}\,(t^{j-i+1};q)_{\la_i-\la_{j+1}}}
{(t^{j-i+1};q)_{\la_i-\la_j}\,(qt^{j-i};q)_{\la_i-\la_{j+1}}}\nonumber\\
\phantom{b_{\la}(q,t)}{}=\prod_{i=1}^n\frac{(t^{n+1-i};q)_{\la_i}}{(qt^{n-i};q)_{\la_i}}\,
\prod_{1\le i<j\le n} \frac{(qt^{j-i},t^{j-i};q)_{\la_i-\la_j}}
{(t^{j-i+1},qt^{j-i-1};q)_{\la_i-\la_j}},\label{bqt}
\end{gather}
for any $n\ge l(\la)$. (It is easy to check that the latter
expression indeed does not depend on $n$.)

The Macdonald polynomials factor nicely under ``principal
specialization'' \cite[p.~343, Example~5]{Ma},
\begin{gather}\label{mcdps}
P_{\la}(1,t,\ldots,t^{N-1};q,t)= t^{\sum_i (i-1)\la_i} \prod_{1\le
i<j\le N} \frac{(t^{j-i+1};q)_{\la_i-\la_j}}
{(t^{j-i};q)_{\la_i-\la_j}}.
\end{gather}
A similar formula holds for principally specialized $Q_{\la}$, by
combining \eqref{mcdnf}, \eqref{bqt}, and \eqref{mcdps}.

We mention two particular useful facts that hold in the case of a
f\/inite set of variables $X=\{x_1,\dots,x_n\}$ (see \cite[p.~323,
Eq.~(4.10), and p.~325, Eq.~(4.17)]{Ma}):
\begin{gather}\label{stab}
P_\la(x_1,\dots,x_n;q,t)=0,\qquad\text{if} \quad l(\la)>n,
\end{gather}
and
\begin{gather}
P_{(\la_1,\la_2,\dots,\la_n)}(x_1,\dots,x_n;q,t)=(x_1\dots
x_n)^{\la_n}\!
P_{(\la_1-\la_n,\la_2-\la_n,\dots,\la_{n-1}-\la_n,0)}(x_1,\dots,x_n;q,t).\!\!\label{fmacd}
\end{gather}

There exists \cite[p.~327]{Ma} an automorphism
$\omega_{q,t}={\omega_{t,q}}^{-1}$ of $\mathsf{Sym}$ such that
\begin{gather}\label{omeg}
\omega_{q,t}(Q_{\la}(q,t))=P_{\la'}(t,q),\qquad
\omega_{q,t}(g_{k}(q,t))=e_{k},
\end{gather}
with $\la'$ the partition conjugate to $\la$, whose parts are
given by $m_k(\la')=\la_k-\la_{k+1}$. In particular \cite[p.~329,
Eq.~(5.5)]{Ma}, the Macdonald polynomials associated with a row or
a column partition are given by
\begin{alignat}{3}
& P_{1^k}(q,t)= e_{k},\qquad &&
Q_{1^k}(q,t)= \frac{(t;t)_{k}}{(q;t)_{k}} \,e_{k},& \notag\\
& P_{(k)}(q,t)= \frac{(q;q)_{k}}{(t;q)_{k}} \, g_{k}(q,t), \qquad
&& Q_{(k)}(q,t)= g_{k}(q,t).& \label{orc}
\end{alignat}
The parameters $q$, $t$ being kept f\/ixed, we shall often write
$P_{\mu}$ or $Q_{\mu}$ for $P_{\mu}(q,t)$ or $Q_{\mu}(q,t)$.

\subsection{Pieri formula}\label{secpieri}

Let $u_1,\ldots,u_n$ be $n$ indeterminates and $\mathbb{N}$ the
set of nonnegative integers. For $\ta =(\ta_1,\ldots,\ta_n)\in
\mathbb{N}^{n}$ let $|\ta|=\sum\limits_{i=1}^{n} \ta_i$ and
def\/ine
\begin{gather*}
d_{\ta_1,\ldots,\ta_n}^{(q,t)} (u_1,\ldots,u_n) =\prod_{k=1}^n
\frac{{(t,q^{|\ta|+1}u_k;q)}_{\ta_k}}
{{(q,q^{|\ta|}tu_k;q)}_{\ta_k}} \prod_{1\le i < j \le n}
\frac{{(tu_i/u_j,q^{-\ta_j+1}u_i/tu_j;q)}_{\ta_i}}
{{(qu_i/u_j,q^{-\ta_j}u_i/u_j;q)}_{\ta_i}}.
\end{gather*}

The Macdonald polynomials satisfy a Pieri formula which
generalizes the classical Pieri formula for Schur
functions~\cite[p.~73, Eq.~(5.16)]{Ma}. This generalization was
obtained by Macdo\-nald~\cite[p.~331]{Ma}, and independently by
Koornwinder~\cite{Ko}.

Most of the time this Pieri formula is stated in combinatorial
terms (as a sum over ``horizontal strips''). In \cite[Th.~3.1]{LS}
it was formulated in ``analytic'' terms:

\begin{Theorem}\label{theopieri}
Let  $\la=(\la_1,\dots,\la_n)$ be an arbitrary partition with
length $\le n$ and $\la_{n+1} \in \mathbb{N}$. Let
$u_i=q^{\la_i-\la_{n+1}}t^{n-i}$, for $1 \le i \le n$. We have
\begin{gather*}
Q_{(\la_1,\ldots,\la_n)} \: Q_{(\la_{n+1})}= \sum_{\ta\in
\mathbb{N}^n} d_{\ta_1,\ldots,\ta_n}^{(q,t)} (u_1,\ldots,u_n)\:
Q_{(\la_1+\ta_1,\ldots,\la_n+\ta_n,\la_{n+1}-|\ta|)}.
\end{gather*}
\end{Theorem}

The Pieri formula def\/ines an inf\/inite transition matrix.
Indeed, let $\mathsf{Sym}(n+1)$ denote the algebra of symmetric
polynomials in $n+1$ independent variables with coef\/f\/icients
in $\mathbb{Q}(q,t)$. Then \cite[p.~313]{Ma} the Macdonald
polynomials $\{Q_\la, \,l(\la) \le n+1\}$ form a basis of
$\mathsf{Sym}(n+1)$, and so do the products $\{Q_{\mu}Q_{(m)},
\,l(\mu) \le n, \,m \ge 0\}$.

\subsection{A recursion formula}\label{secmain}

Again, let $u_1,\ldots,u_n$ be $n$ indeterminates. For $\ta
=(\ta_1,\ldots,\ta_n)\in \mathbb{N}^{n}$ def\/ine
\begin{gather*}
c^{(q,t)}_{\ta_1,\ldots,\ta_n} (u_1,\dots,u_n)={} \prod_{i=1}^n
t^{\ta_i} \,\frac{(q/t,qu_k;q)_{\ta_i}}{(q,qtu_i;q)_{\ta_i}}\,
\prod_{1\le i < j \le n}
\frac{{(qu_i/tu_j,q^{-\ta_j}tu_i/u_j;q)}_{\ta_i}}
{{(qu_i/u_j,q^{-\ta_j}u_i/u_j;q)}_{\ta_i}}\\
\qquad{}\times\frac{1}{\prod\limits_{1\le i < j \le
n}(q^{\ta_i}u_i-q^{\ta_j}u_j)} \, \det_{1\le i,j \le n}\!
\left[\big(q^{\ta_i}u_i\big)^{n-j} \left(1-t^{j-1}
\frac{1-q^{\ta_i}tu_i}{1-q^{\ta_i}u_i} \prod_{s=1}^n
\frac{u_s-q^{\ta_i}u_i}{tu_s-q^{\ta_i}u_i}\right)\right].
\end{gather*}

The following recursion formula for Macdonald polynomials was
proved in \cite[Th.~4.1]{LS} by inverting the Pieri formula in
Theorem~\ref{theopieri} utilizing a special case of the
multidimensional matrix inverse in Corollary~\ref{cormi2}.

\begin{Theorem}\label{theomain}
Let  $\la=(\la_1,\dots ,\la_{n+1})$ be an arbitrary partition with
length $\le n+1$. Let $u_i=q^{\la_i-\la_{n+1}}t^{n-i}$, for $1\le
i \le n$. We have
\begin{gather*}
Q_{(\la_1,\ldots,\la_{n+1})}= \sum_{\ta\in\mathbb{N}^n}
c^{(q,t)}_{\ta_1,\ldots,\ta_n} (u_1,\ldots,u_n)\:
Q_{(\la_{n+1}-|\ta|)} \: Q_{(\la_1+\ta_1,\ldots,\la_n+\ta_n)}.
\end{gather*}
\end{Theorem}

In the case $n=1$, i.e.\ for partitions of length two,
Theorem~\ref{theomain} reads
\begin{gather*}
Q_{(\la_1,\la_2)}= \sum_{\ta\in\mathbb{N}} c^{(q,t)}_{\ta} (u) \,
Q_{(\la_{2}-\ta)} \: Q_{(\la_1+\ta)},
\end{gather*}
with $u=q^{\la_1-\la_2}$ and
\begin{gather*}
c^{(q,t)}_{\ta} (u)= t^{\ta}
\,\frac{(q/t,qu;q)_{\ta}}{(q,qtu;q)_{\ta}}\, \left(1-
\frac{(1-tq^{\ta}u)}{(1-q^{\ta}u)}
\frac{(u-q^{\ta}u)}{(tu-q^{\ta}u)}\right)\\
\phantom{c^{(q,t)}_{\ta} (u)}{}= t^{\ta}
\,\frac{(q/t,qu;q)_{\ta}}{(q,qtu;q)_{\ta}}\,
\frac{(t-1)}{(t-q^{\ta})}\frac{(1-q^{2\ta}u)}{(1-q^{\ta}u)} =
t^{\ta} \,\frac{(1/t,u;q)_{\ta}}{(q,qtu;q)_{\ta}}\,
\frac{(1-q^{2\ta}u)}{(1-u)}.
\end{gather*}
This special case is due to Jing and J\'ozef\/iak~\cite{JJ}. On
the other hand, in the case $n=2$, i.e.\ for partitions of length
three, Theorem~\ref{theomain} reduces to a formula stated by
Lassalle~\cite{La}.

Application of the automorphism $\omega_{q,t}$ to
Theorem~\ref{theomain}, while taking into account \eqref{omeg},
gives the following equivalent result (cf.\ \cite[Th.~4.2]{LS}).
\begin{Theorem}\label{theodual}
Let  $\la=(1^{m_1},2^{m_2},\ldots,(n+1)^{m_{n+1}})$ be an
arbitrary partition consisting of parts at most equal to $n+1$.
Def\/ine $u_i=q^{n-i}t^{\sum\limits_{j=i}^n m_j}$, for $1\le i \le
n+1$ . We have
\begin{gather*}
P_\la= \sum_{\ta\in\mathbb{N}^n} c^{(t,q)}_{\ta_1,\ldots,\ta_n}
(u_1,\ldots,u_n) \: e_{m_{n+1}-|\ta|} \:
P_{(1^{m_1+\ta_1-\ta_2},\ldots,
{(n-1)}^{m_{n-1}+\ta_{n-1}-\ta_n},n^{m_n+m_{n+1}+\ta_n})}.
\end{gather*}
\end{Theorem}

\section{Basic hypergeometric series}\label{secbhs}

\subsection{Classical (one-dimensional) basic hypergeometric series}

First we recall some standard notations for $q$-series and basic
hypergeometric series (see Gasper and Rahman's text~\cite{GR}, for
a standard reference). In the following we shall consider $q$ a
(f\/ixed) complex parameter with $0<|q|<1$. Further, for a complex
parameter $a$, we use the def\/inition of the $q$-shifted
factorial as given in \eqref{ipr} and \eqref{cdn}, for any complex
number $k$.

The basic hypergeometric ${}_r\phi_{r-1}$ series with upper
parameters $a_1,\dots,a_r$, lower parameters $b_1,\dots,b_{r-1}$,
base $q$, and argument $z$ is def\/ined as follows:

\begin{gather*}
{}_r\phi_{r-1}\!\left[\begin{matrix}a_1,a_2,\dots,a_r\\
b_1,b_2,\dots,b_{r-1}\end{matrix}\,;q,z\right]= \sum _{k=0}
^{\infty}\frac {(a_1,a_2,\dots,a_r;q)_k}
{(q,b_1,\dots,b_{r-1};q)_k}\,z^k.
\end{gather*}

An ${}_r\phi_{r-1}$ series terminates if one of its upper
parameters is of the form $q^{-m}$ with $m=0,1,2,\dots$, because
\begin{gather*}
(q^{-m};q)_k=0,\qquad k=m+1,m+2,\dots.
\end{gather*}
The ${}_r\phi_{r-1}$ series, if it does not terminate, converges
absolutely for $|z|<1$.

In our computations throughout this paper, we make frequent use of
some elementary identities for $q$-shifted factorials, listed in
\cite[Appendix~I]{GR}. We list a few important summation theorems
from \cite[Appendix~II]{GR}, for quick reference.

We start with the (nonterminating) $q$-binomial theorem (cf.\
\cite[Eq.~(II.2)]{GR})
\begin{gather}\label{sum10}
{}_1\phi_0\!\left[\begin{matrix}a\\-\end{matrix};q,z\right]=
\frac{(az;q)_\infty}{(z;q)_\infty},\qquad |z|<1,
\end{gather}
an identity f\/irst discovered by Cauchy~\cite{C} in 1843.

Further, we have the $q$-Gau{\ss} summation (cf.\
\cite[Eq.~(II.8)]{GR}),
\begin{gather}\label{sum21}
{}_2\phi_1\!\left[\begin{matrix}a,b\\c\end{matrix};q,\frac
c{ab}\right]= \frac{(c/a,c/b;q)_\infty}{(c,c/ab;q)_\infty},\qquad
|c/ab|<1,
\end{gather}
which is due to Heine~\cite{H}.

An identity of fundamental importance is the terminating
${}_6\phi_5$ summation theorem (cf.\ \cite[Eq.~(II.21)]{GR})
\begin{gather}\label{sum65}
{}_6\phi_5\!\left[\begin{matrix}a,\,q\sqrt{a},-q\sqrt{a},b,c,q^{-m}\\
\sqrt{a},-\sqrt{a},aq/b,aq/c,aq^{1+m}\end{matrix}\,;
q,\frac{aq^{1+m}}{bc}\right]=
\frac{(aq,aq/bc;q)_m}{(aq/b,aq/c;q)_m}.
\end{gather}
This identity can be extended (by analytic continuation) to
Rogers'~\cite{Ro} nonterminating ${}_6\phi_5$ summation theorem
(cf.\ \cite[Eq.~(II.20)]{GR})
\begin{gather}\label{sumnt65}
{}_6\phi_5\!\left[\begin{matrix}a,\,q\sqrt{a},-q\sqrt{a},b,c,d\\
\sqrt{a},-\sqrt{a},aq/b,aq/c,aq/d\end{matrix}\,;
q,\frac{aq}{bcd}\right]= \frac{(aq,aq/bc,aq/bd,aq/cd;q)_\infty}
{(aq/b,aq/c,aq/d,aq/bcd;q)_\infty},
\end{gather}
where $|aq/bcd|<1$. Clearly, \eqref{sumnt65} reduces to
\eqref{sum65} for $d=q^{-m}$.

The series in \eqref{sum65} and \eqref{sumnt65} both are examples
of so-called ``very-well-poised'' series. Another example is
Jackson's~\cite{Ja} terminating very-well-poised ${}_8\phi_7$
summation (cf.\ \cite[Eq.~(II.22)]{GR}),
\begin{gather}\label{sum87}
{}_8\phi_7\!\!\left[\begin{matrix}a,\,q\sqrt{a},-q\sqrt{a},b,c,d,
a^2q^{1+n}/bcd,q^{-n}\\
\sqrt{a},-\sqrt{a},aq/b,aq/c,aq/d,bcdq^{-n}/a,aq^{1+n}\end{matrix}\,;
q,q\right]\!= \frac{(aq,aq/bc,aq/bd,aq/cd;q)_n}
{(aq/b,aq/c,aq/d,aq/bcd;q)_n},\!\!\!
\end{gather}
which contains \eqref{sumnt65} as a special case ($n\to\infty$).
(The series in \eqref{sum87} is even ``balanced'' but we do not
need this notion here.) Jackson's summation \eqref{sum87} stands
on the top of the classical hierarchy of summation theorems for
very-well-poised basic hypergeometric series.

An ${}_r\phi_{r-1}$ series is well-poised if the parameters
satisfy the relations
\begin{gather*}
qa_1=a_2b_1=a_3b_2=\dots=a_rb_{r-1}.
\end{gather*}
It is very-well-poised if, in addition, $a_2=q\sqrt{a_1}$,
$a_3=-q\sqrt{a_1}$. Note that
\begin{gather}\label{vwp}
\frac{(q\sqrt{a_1},-q\sqrt{a_1};q)_k}
{(\sqrt{a_1},-\sqrt{a_1};q)_k} =\frac{1-a_1q^{2k}}{1-a_1}
\end{gather}
appears as a factor in the summand of a very-well-poised series.
The parameter $a_1$ is usually referred to as the special
parameter of such a series, and we call \eqref{vwp} the
very-well-poised term of the series.

\subsection[$A_{n-1}$ basic hypergeometric series]{$\boldsymbol{A_{n-1}}$ basic hypergeometric series} \label{secanseries}

For convenience, we employ here (and elsewhere) the notation
$|{\mathbf k}|= k_1+\dots+k_n$, where ${\mathbf
k}=(k_1,\dots,k_n)$. We recall the following fundamental result
from \cite[Th.~1.49]{Mi1}.
\begin{Proposition}[Milne]\label{ftans}
Let $a_1,\dots,a_n$, and $u_1,\dots,u_n$ be indeterminate, let $M$
be a nonnegative integer. Then
\begin{gather}\label{fteq}
\sum_{\underset{|{\mathbf k}|=M}{k_1,\dots,k_n\ge 0}} \prod_{1\le
i<j\le n}\frac {u_iq^{k_i}-u_jq^{k_j}} {u_i-u_j}
\prod_{i,j=1}^n\frac {(a_ju_i/u_j;q)_{k_i}} {(qu_i/u_j;q)_{k_i}}=
\frac {(a_1a_2\dots a_n;q)_M}{(q;q)_M}.
\end{gather}
\end{Proposition}
The $n=2$ case of \eqref{fteq} is equivalent to \eqref{sum65}.
Proposition~\ref{ftans} played an important role in Milne's
elementary derivation of the Macdonald identities for the
af\/f\/ine Lie algebra $A_n^{(1)}$. It further serves as a
starting point in Milne's approach to $A_{n-1}$ or $U(n)$ series,
see \cite{Mi6}. Its original proof made use of a particular
$q$-dif\/ference equation, and induction. An even simpler
inductive proof of Proposition~\ref{ftans} was recently discovered
by Rosengren~\cite[Sec.~2]{R1} (see also \cite[Sec.~11.7]{GR}).

An immediate consequence of Proposition~\ref{ftans} is the
following $n$-dimensional extension of \eqref{sum10}, f\/irst
given in \cite[Th.~1.47]{Mi1}.

\begin{Proposition}[(Milne) An
$A_{n-1}$ nonterminating $q$-binomial theorem]\label{an10} Let
$a_1,\dots,a_n$, $z$, and $u_1,\dots,u_n$ be indeterminate. Then
\begin{gather*}
\sum_{k_1,\dots,k_n=0}^{\infty} \prod_{1\le i<j\le n}\frac
{u_iq^{k_i}-u_jq^{k_j}} {u_i-u_j} \prod_{i,j=1}^n\frac
{(a_ju_i/u_j;q)_{k_i}} {(qu_i/u_j;q)_{k_i}}\cdot z^{|{\mathbf
k}|}= \frac {(a_1\dots a_nz;q)_{\infty}}{(z;q)_{\infty}},
\end{gather*}
provided $|z|<1$.
\end{Proposition}
Clearly, Proposition~\ref{an10} follows from
Proposition~\ref{ftans} by multiplying both sides of \eqref{fteq}
by~$z^M$ and summing over all $M\ge 0$. The right-hand side is
then simplif\/ied using \eqref{sum10}.

In the above cases, we have
\begin{gather}\label{anvandy}
\prod_{1\le i<j\le n}\frac {u_iq^{k_i}-u_jq^{k_j}} {u_i-u_j}
\end{gather}
appearing as a factor of the summand. Since we may associate
\eqref{anvandy} with the product side of the Weyl denominator
formula for the root system $A_{n-1}$ (see e.g.\
D.~Stanton~\cite{St}), we call such a series $A_{n-1}$ basic
hypergeometric series. Note that often in the literature these
$n$-dimensional series are called $A_n$ series instead of
$A_{n-1}$ series.

Note that on the left-hand side of \eqref{an65eq}, e.g., we have,
in addition to \eqref{anvandy},
\begin{gather}\label{mvwp}
\prod_{i=1}^n\frac{1-au_iq^{k_i+|{\mathbf k}|}}{1-au_i}
\end{gather}
appearing as a factor in the summand of the series. It is easy to
see that the $n=1$ case of \eqref{mvwp} essentially reduces to
\eqref{vwp}. We call the product of \eqref{anvandy} and
\eqref{mvwp} the very-well-poised term of the multiple series. To
clarify the special appearance of the very-well-poised term (even
in the one-dimensional case), it is useful to view the series in
one higher dimension. In particular, we can write
\begin{gather*}
\prod_{1\le i<j\le n}\frac {u_iq^{k_i}-u_jq^{k_j}}{u_i-u_j}
\prod_{i=1}^n\frac{1-au_iq^{k_i+|{\mathbf k}|}}{1-au_i}
=q^{|{\mathbf k}|} \prod_{1\le i<j\le n+1}\frac
{u_iq^{k_i}-u_jq^{k_j}}{u_i-u_j},
\end{gather*}
where $u_{n+1}=1/a$ and $k_{n+1}=-|{\mathbf k}|$. Indeed, some
$A_{n-1}$ basic hypergeometric series identities are sometimes
better viewed as identities associated to the {\em affine} root
system $\tilde{A}_n$ (or, equivalently, the special unitary group
$SU(n+1)$).

Multidimensional basic hypergeometric series associated with the
root system $A_{n-1}$ (or equi\-valently, associated with the
unitary group $U(n)$) have their origin in the work of the three
mathe\-matical physicists Holman, Biedenharn, and Louck, see
\cite{Ho,HBL}. Their work was done in the context of the quantum
theory of angular momentum, using methods relying on the
representation theory of the unitary group $U(n)$. Subsequently,
extensive investigations in the theory of multiple basic
hypergeometric series associated with the root system $A_{n-1}$
have been carried out by R.A.~Gustafson, S.C.~Milne, and several
other researchers. As result, many of the classical formulae for
basic hypergeometric series from \cite{GR} have already been
generalized to the setting of $A_{n-1}$ series, see e.g.\ the
references \cite{BM,DG,G1,G2,Mi1,Mi2,Mi3,Mi4,Mi5,ML,MN,R2}.

An immediate consequence of Proposition~\ref{ftans} is the
following extension of \eqref{sum65}, a result f\/irst obtained in
\cite[Th.~1.38]{Mi2}.

\begin{Proposition}[(Milne) An
$A_{n-1}$ terminating ${}_6\phi_5$ summation]\label{an65} Let $a$,
$b$, $c_1,\dots,c_n$, and $u_1,\dots,u_n$ be indeterminate, let
$M$ be a nonnegative integer. Then
\begin{gather}
\sum_{\underset{|{\mathbf k}|\le M}{k_1,\dots,k_n\ge 0}}
\prod_{1\le i<j\le n}\frac {u_iq^{k_i}-u_jq^{k_j}} {u_i-u_j}
\prod_{i=1}^n\frac{1-au_iq^{k_i+|{\mathbf k}|}}{1-au_i}
\prod_{i=1}^n\frac {(au_i;q)_{|{\mathbf k}|}}
{(au_iq/c_i;q)_{|{\mathbf k}|}}\nonumber\\
\qquad{}\times \prod_{i,j=1}^n\frac {(c_ju_i/u_j;q)_{k_i}}
{(qu_i/u_j;q)_{k_i}} \prod_{i=1}^n\frac {(bu_i;q)_{k_i}}
{(au_iq^{1+M};q)_{k_i}}\cdot \frac{(q^{-M};q)_{|{\mathbf
k}|}}{(aq/b;q)_{|{\mathbf k}|}}
\left(\frac{aq^{1+M}}{bc_1\cdots c_n}\right)^{|{\mathbf k}|}\nonumber\\
\qquad{}=\frac{(aq/bc_1\cdots c_n;q)_M}{(aq/b;q)_M}
\prod_{i=1}^n\frac{(au_iq;q)_M}{(au_iq/c_i;q)_M}.\label{an65eq}
\end{gather}
\end{Proposition}

\begin{proof}[Proof (see e.g.\ \cite{Bh})]
To derive Proposition~\ref{an65} from Proposition~\ref{ftans},
just replace $n$ by $n+1$ in \eqref{fteq}, and replace $k_{n+1}$
by $M-|{\mathbf k}|$, so that in the resulting sum the summation
index ranges over the $n$-tetrahedron $0\le|{\mathbf k}|\le M$.
Finally, suitably relabel parameters to obtain \eqref{an65eq}.
\end{proof}

Proposition~\ref{an65} can be extended by analytic continuation to
the following multivariable gene\-ralization of \eqref{sumnt65}, a
result f\/irst obtained in \cite[Th.~1.44]{Mi2}.

\begin{Proposition}[(Milne) An
$A_{n-1}$ nonterminating ${}_6\phi_5$ summation]\label{annt65} Let
$a$, $b$, $c_1,\dots,c_n$, $d$, and $u_1,\dots,u_n$ be
indeterminate. Then
\begin{gather}
\sum_{k_1,\dots,k_n\ge 0} \prod_{1\le i<j\le n}\frac
{u_iq^{k_i}-u_jq^{k_j}} {u_i-u_j}
\prod_{i=1}^n\frac{1-au_iq^{k_i+|{\mathbf k}|}}{1-au_i}
\prod_{i=1}^n\frac {(au_i;q)_{|{\mathbf k}|}}
{(au_iq/c_i;q)_{|{\mathbf k}|}}\nonumber\\
\qquad{}\times \prod_{i,j=1}^n\frac {(c_ju_i/u_j;q)_{k_i}}
{(qu_i/u_j;q)_{k_i}} \prod_{i=1}^n\frac {(bu_i;q)_{k_i}}
{(au_iq/d;q)_{k_i}}\cdot \frac{(d;q)_{|{\mathbf
k}|}}{(aq/b;q)_{|{\mathbf k}|}}
\left(\frac{aq}{bc_1\cdots c_nd}\right)^{|{\mathbf k}|}\nonumber\\
\qquad{}=\frac{(aq/bc_1\cdots c_n,aq/bd;q)_\infty}
{(aq/b,aq/bc_1\cdots c_nd;q)_\infty}
\prod_{i=1}^n\frac{(au_iq,au_iq/c_id;q)_\infty}
{(au_iq/c_i,au_iq/d;q)_\infty},\label{annt65eq}
\end{gather}
provided $|aq/bc_1\cdots c_nd|<1$.
\end{Proposition}
Clearly, Proposition~\ref{annt65} reduces to
Proposition~\ref{an65} for $d=q^{-M}$.
\begin{proof}
Both sides of \eqref{annt65eq} are analytic in $1/d$ in a domain
$\mathcal D$ around the origin. Now, if $1/d=q^M$ the formula
holds due to Proposition~\ref{an65}, for all $M=0,1,2,\dots$.
Since the inf\/inite sequence $q^M$ has an accumulation point,
namely $0$, in the domain of analyticity of $1/d$, we can apply
the identity theorem to deduce that \eqref{annt65eq} is true for
all $1/d$ throughout $\mathcal D$. By analytic continuation the
identity holds for all $1/d$ in its region of convergence, i.e.\
in the disc $|1/d|<|bc_1\cdots c_n/aq|$.
\end{proof}

We will need the following identity from \cite[Th.~7.6]{Mi5} which
can be obtained from Proposition~\ref{annt65} by letting $b\mapsto
aq/c$, followed by $a\to 0$, and the subsequent substitutions
$c_i\mapsto a_i$, $i=1,\dots,n$, and $d\mapsto b$. It represents a
multivariable generalization of the $q$-Gau{\ss} summation in
\eqref{sum21}.
\begin{Proposition}[(Milne) An
$A_{n-1}$ nonterminating ${}_2\phi_1$ summation]\label{annt21} Let
$a_1,\dots,a_n$, $b$, $c$, and $u_1,\dots,u_n$ be indeterminate.
Then
\begin{gather*}
\sum_{k_1,\dots,k_n\ge 0} \prod_{1\le i<j\le n}\frac
{u_iq^{k_i}-u_jq^{k_j}}{u_i-u_j} \prod_{i,j=1}^n\frac
{(a_ju_i/u_j;q)_{k_i}} {(qu_i/u_j;q)_{k_i}}\cdot
\frac{(b;q)_{|{\mathbf k}|}}{(c;q)_{|{\mathbf k}|}}
\left(\frac{c}{a_1\cdots a_nb}\right)^{|{\mathbf k}|}\\
\qquad{}=\frac{(c/a_1\cdots a_n,c/b;q)_\infty} {(c,c/a_1\cdots
a_nb;q)_\infty},
\end{gather*}
provided $|c/a_1\cdots a_nb|<1$.
\end{Proposition}

The following identity in Proposition~\ref{an87} is even more
general than Propositions~\ref{an65} and \ref{annt65}. It was
f\/irst obtained in \cite[Th.~6.17]{Mi7}. A simple inductive proof
of Proposition~\ref{an87} was recently discovered by
Rosengren~\cite[Sec.~5]{R1} (see also~\cite[Sec.~11.7]{GR}).

\begin{Proposition}[(Milne) An
$A_{n-1}$ terminating ${}_8\phi_7$ summation]\label{an87} Let $a$,
$b$, $c_1,\dots,c_n$, $d$, and $u_1,\dots,u_n$ be indeterminate,
let $M$ be a nonnegative integer. Then
\begin{gather*}
\sum_{\underset{|{\mathbf k}|\le M}{k_1,\dots,k_n\ge 0}}
\prod_{1\le i<j\le n}\frac {u_iq^{k_i}-u_jq^{k_j}}{u_i-u_j}
\prod_{i=1}^n\frac{1-au_iq^{k_i+|{\mathbf k}|}}{1-au_i}
\prod_{i=1}^n\frac {(au_i;q)_{|{\mathbf k}|}}
{(au_iq/c_i;q)_{|{\mathbf k}|}} \prod_{i,j=1}^n\frac
{(c_ju_i/u_j;q)_{k_i}}
{(qu_i/u_j;q)_{k_i}}\nonumber\\
\qquad{}\times \prod_{i=1}^n\frac {(bu_i,a^2u_iq^{1+M}/bc_1\cdots
c_nd;q)_{k_i}} {(au_iq/d,au_iq^{1+M};q)_{k_i}}\cdot
\frac{(d,q^{-M};q)_{|{\mathbf k}|}}
{(aq/b,bc_1\cdots c_ndq^{-M}/a;q)_{|{\mathbf k}|}}\,q^{|{\mathbf k}|}\nonumber\\
\qquad{}=\frac{(aq/bd,aq/bc_1\cdots c_n;q)_M}{(aq/b,aq/bc_1\cdots
c_nd;q)_M}\, \prod_{i=1}^n\frac{(au_iq,au_iq/c_id;q)_M}
{(au_iq/c_i,au_iq/d;q)_M}.
\end{gather*}
\end{Proposition}

For convenience, we also list the following identity which is not
of $A_{n-1}$ but of $C_n$ type. It was derived in
\cite[Th.~4.1]{DG}, and, independently, in \cite[Th.~6.13]{ML}.

\begin{Proposition}[(Denis--Gustafson; Milne--Lilly)
A $C_n$ terminating ${}_8\phi_7$ summation]\label{cn87} Let $a$,
$b$, $c$, $d$ and $u_1,\dots,u_n$ be indeterminate, let
$m_1,\dots,m_n$ be nonnegative integers. Then
\begin{gather*}
\underset{i=1,\dots,n}{\sum_{0\le k_i\le m_i}} \prod_{1\le i<j\le
n}\frac {u_iq^{k_i}-u_jq^{k_j}}{u_i-u_j} \prod_{1\le i\le j\le
n}\frac{1-au_iu_jq^{k_i+k_j}}{1-au_iu_j}
\prod_{i,j=1}^n\frac{(q^{-m_j}u_i/u_j,au_iu_j;q)_{k_i}}
{(au_iu_jq^{1+m_j},qu_i/u_j;q)_{k_i}}\nonumber\\
\qquad{}\times
\prod_{i=1}^n\frac{(bu_i,cu_i,du_i,a^2u_iq^{1+|\mathbf
m|}/bcd;q)_{k_i}} {(au_iq/b,au_iq/c,au_iq/d,bcdu_iq^{-|\mathbf
m|}/a;q)_{k_i}}\cdot
q^{|{\mathbf k}|}\nonumber\\
\qquad{}=\prod_{1\le i<j\le n}(au_iu_jq;q)_{m_i+m_j}^{-1}
\prod_{i,j=1}^n(au_iu_jq;q)_{m_i}\nonumber\\
\qquad{}\times \frac{(aq/bc,aq/bd,aq/cd;q)_{|\mathbf m|}}
{\prod\limits_{i=1}^n(au_iq/b,au_iq/c,au_iq/d,aq^{1+|\mathbf m|-m_i}/bcdu_i;q)_{m_i}}.
\end{gather*}
\end{Proposition}

\section{Multivariable basic hypergeometric series
of a new type}\label{secmsnt}

Here we derive several new summations for multivariable basic
hypergeometric series of a new type. Their special feature is that
the series involve specif\/ic determinants in the summand. The
f\/irst result is obtained from the $A_{n-1}$ ${}_8\phi_7$
summation in Proposition~\ref{an87}, combined with the
multivariable matrix inverse in Corollary~\ref{cormi2}.

\begin{Theorem}[A multivariable
terminating ${}_8\phi_7$ summation]\label{an87n} Let $b$, $d$,
$t_0,t_1,\dots,t_n$, $u_1,\dots,u_n$ be indeterminate, let $M$ be
a nonnegative integer. Then
\begin{gather}
\sum_{\underset{|{\mathbf k}|\le M}{k_1,\dots,k_n\ge 0}}
\prod_{i,j=1}^n\frac{{(qu_i/t_iu_j;q)}_{k_i}}
{{(qu_i/u_j;q)}_{k_i}}\, \prod_{1\le i<j\le
n}\frac{(t_ju_i/u_j;q)_{k_i-k_j}} {(qu_i/t_iu_j;q)_{k_i-k_j}}
\prod_{1\le i < j \le n}\frac 1{u_i-u_j}\nonumber\\
\qquad{}\times \det_{1\le i,j \le n}\! \left[(u_iq^{k_i})^{n-j}
\left(1-t_i^{j-n-1}\frac{1-t_0u_iq^{k_i}} {1-t_0u_iq^{k_i}/t_i}
\prod_{s=1}^n \frac{u_iq^{k_i}-u_s}
{u_iq^{k_i}/t_i-u_s}\right)\right]\nonumber\\
\qquad{}\times \frac{(d,q^{-M};q)_{|\mathbf k|}}
{(bdq^{-M}/t_0,t_0q/bt_1\cdots t_n;q)_{|\mathbf k|}}\,
\prod_{i=1}^n \frac{(t_0u_iq/t_i,bu_i,t_0^2u_iq^{1+M}/bdt_1\cdots
t_n;q)_{k_i}}
{(t_0u_iq,t_0u_iq/dt_i,t_0u_iq^{1+M}/t_i;q)_{k_i}}\nonumber\\
\qquad{}\times \prod_{i=1}^n \frac{(dq^{-M}/t_0u_i;q)_{|\mathbf
k|-k_i}} {(dt_iq^{-M}/t_0u_i;q)_{|\mathbf k|-k_i}}\cdot
q^{|\mathbf k|+\sum\limits_{i=1}^n(1-i)k_i}\,
\prod_{i=1}^nt_i^{(i-1)k_i+\sum\limits_{j=i+1}^nk_j}\nonumber\\
\qquad{}= \frac{(t_0q/b,t_0q/bdt_1\cdots t_n;q)_M}
{(t_0q/bd,t_0q/bt_1\cdots t_n;q)_M}
\prod_{i=1}^n\frac{(t_0u_iq/t_i,t_0u_iq/d;q)_M}
{(t_0u_iq,t_0u_iq/dt_i;q)_M}.\label{an87neq}
\end{gather}
\end{Theorem}

\begin{proof}
We apply multidimensional inverse relations. We have
\eqref{rotinv2} by the $a\mapsto t_0q^{|\mathbf l|}$, $c_i\mapsto
t_i$, $d\mapsto dq^{|\mathbf l|}$, $u_i\mapsto u_iq^{l_i}$,
$i=1,\dots,n$, and $M\mapsto M-|\mathbf l|$ case of
Proposition~\ref{an87}, where
\begin{gather*}
b_{\mathbf k}=\prod_{1\le i < j \le n}(u_iq^{k_i}-u_jq^{k_j})\,
\prod_{i,j=1}^n\frac{{(t_ju_i/u_j;q)}_{k_i-k_j}}
{{(qu_i/u_j;q)}_{k_i-k_j}}\,\prod_{i=1}^n\frac{(t_0u_iq;q)_{k_i+|\mathbf
k|}}
{(t_0u_iq/t_i;q)_{k_i+|\mathbf k|}}\\
\phantom{b_{\mathbf k}=}{}\times
\prod_{i=1}^n\frac{(bu_i,t_0^2u_iq^{1+M}/bdt_1\cdots t_n;q)_{k_i}}
{(t_0u_iq/d,t_0u_iq^{1+M};q)_{k_i}}\cdot
\frac{(d,q^{-M};q)_{|\mathbf k|}} {(t_0q/b,bdt_1\cdots
t_nq^{-M}/t_0;q)_{|\mathbf k|}}\,q^{|\mathbf k|}
\end{gather*}
and
\begin{gather*}
a_{\mathbf l}=\frac{(t_0q/bd,t_0q/bt_1\cdots t_n;q)_M}
{(t_0q/b,t_0q/bdt_1\cdots t_n;q)_M}\,\prod_{i=1}^n
\frac{(t_0u_iq,t_0u_iq/dt_i;q)_M}{(t_0u_iq/t_i,t_0u_iq/d;q)_M}\\
\phantom{a_{\mathbf l}=}{}\times \prod_{1\le i < j \le
n}(u_iq^{l_i}-u_jq^{l_j})\,
\prod_{i,j=1}^n\frac{{(t_ju_i/u_j;q)}_{l_i-l_j}}
{{(qu_i/u_j;q)}_{l_i-l_j}}\,\prod_{i=1}^n
\frac{(dq^{-M}/t_0u_i;q)_{|\mathbf l|-l_i}}
{(dt_iq^{-M}/t_0u_i;q)_{|\mathbf l|-l_i}}\\
\phantom{a_{\mathbf l}=}{}\times
\prod_{i=1}^n\frac{(bu_i,t_0^2u_iq^{1+M}/bdt_1\cdots t_n;q)_{l_i}}
{(t_0u_iq/dt_i,t_0u_iq^{1+M}/t_i;q)_{l_i}}\cdot
\frac{(d,q^{-M};q)_{|\mathbf l|}} {(t_0q/bt_1\cdots
t_n,bdq^{-M}/t_0;q)_{|\mathbf l|}}\, q^{|\mathbf
l|}\,\prod_{i=1}^nt_i^{-l_i},
\end{gather*}
and $g_{\mathbf k\mathbf l}$ as in \eqref{cormig}. Therefore, we
must have \eqref{rotinv1} with the above sequences $a_{\mathbf l}$
and $b_{\mathbf k}$, and $f_{\mathbf m\mathbf k}$ as in
\eqref{cormif}. After the substitutions $t_0\mapsto
t_0q^{-|\mathbf k|}$, $d\mapsto dq^{-|\mathbf k|}$, $u_i\mapsto
u_iq^{-k_i}$, $m_i\mapsto m_i+k_i$, $i=1,\dots,n$, $M\mapsto
M+|\mathbf k|$, simplif\/ications and subsequent relabelling
($m_i\mapsto k_i$, $i=1,\dots,n$), we obtain \eqref{an87neq}.
\end{proof}

The following identity is obtained from Theorem~\ref{an87n} by a
polynomial argument.

\begin{Corollary}[A multivariable
terminating ${}_8\phi_7$ summation]\label{an87np} Let $a$, $b$,
$c$, $d$, and $u_1,\dots,u_n\!$ be indeterminate, let
$m_1,\dots,m_n$ be nonnegative integers. Then
\begin{gather}
\sum_{\underset{i=1,\dots,n}{0\le k_i\le m_i}}
\prod_{i,j=1}^n\frac{{(q^{1-m_i}u_i/u_j;q)}_{k_i}}
{{(qu_i/u_j;q)}_{k_i}}\, \prod_{1\le i<j\le
n}\frac{(q^{m_j}u_i/u_j;q)_{k_i-k_j}}
{(q^{1-m_i}u_i/u_j;q)_{k_i-k_j}}
\prod_{1\le i < j \le n}\frac 1{u_i-u_j}\nonumber\\
\qquad{}\times \det_{1\le i,j \le n}\! \left[(u_iq^{k_i})^{n-j}
\left(1-q^{(j-n-1)m_i}\frac{1-au_iq^{k_i+|\mathbf m|}}
{1-au_iq^{k_i+|\mathbf m|-m_i}} \prod_{s=1}^n
\frac{u_iq^{k_i}-u_s}
{u_iq^{k_i-m_i}-u_s}\right)\right]\nonumber\\
\qquad{}\times \frac{(c,d;q)_{|\mathbf k|}} {(bcdq^{-|\mathbf
m|}/a,aq/b;q)_{|\mathbf k|}}\, \prod_{i=1}^n
\frac{(au_iq^{1+|\mathbf m|-m_i},bu_i,a^2u_iq^{1+|\mathbf
m|}/bcd;q)_{k_i}} {(au_iq^{1+|\mathbf m|},au_iq^{1+|\mathbf
m|-m_i}/c,
au_iq^{1+|\mathbf m|-m_i}/d;q)_{k_i}}\nonumber\\
\qquad{}\times \prod_{i=1}^n \frac{(cdq^{-|\mathbf
m|}/au_i;q)_{|\mathbf k|-k_i}} {(cdq^{m_i-|\mathbf
m|}/au_i;q)_{|\mathbf k|-k_i}}\cdot q^{|\mathbf k|+
\sum_{i=1}^n\big((1-i)k_i+n_i[(i-1)k_i+\sum_{j=i+1}^nk_j]\big)}\nonumber\\
\qquad{}= \frac{(aq/bc,aq/bd;q)_{|\mathbf m|}}
{(aq/b,aq/bcd;q)_{|\mathbf m|}}\,
\prod_{i=1}^n\frac{(au_iq^{1+|\mathbf m|-m_i}, au_iq^{1+|\mathbf
m|-m_i}/cd;q)_{m_i}} {(au_iq^{1+|\mathbf
m|-m_i}/c,au_iq^{1+|\mathbf m|-m_i}/d;q)_{m_i}}.\label{an87npeq}
\end{gather}
\end{Corollary}

\begin{proof}
First we write the right hand-side of \eqref{an87neq} as quotient
of inf\/inite products using \eqref{ipr}. Then by the
$t_0=aq^{|\mathbf m|}$ and $t_i=q^{m_i}$, $i=1,\dots,n$, case of
Theorem~\ref{an87} it follows that the identity \eqref{an87npeq}
holds for $c=q^{-M}$. By clearing out denominators in
\eqref{an87npeq}, we get a polynomial equation in $c$, which is
true for $q^{-M}$, $M=0,1,\dots$, and thus obtain an identity in
$c$.
\end{proof}

The following multivariable nonterminating ${}_6\phi_5$ summation
is obtained by letting $M\to\infty$ in Theorem~\ref{an87n} (while
appealing to Tannery's theorem (cf.\ \cite{Bw}) for
justif\/ication of taking term-wise limits).

\begin{Corollary}[A multivariable
nonterminating ${}_6\phi_5$ summation]\label{an65nt} Let $b$, $d$,
$t_0,t_1,\dots,t_n$, and $u_1,\dots,u_n$ be indeterminate. Then
\begin{gather*}
\sum_{k_1,\dots,k_n\ge 0}
\prod_{i,j=1}^n\frac{{(qu_i/t_iu_j;q)}_{k_i}}
{{(qu_i/u_j;q)}_{k_i}}\, \prod_{1\le i<j\le
n}\frac{(t_ju_i/u_j;q)_{k_i-k_j}} {(qu_i/t_iu_j;q)_{k_i-k_j}}
\prod_{1\le i < j \le n}\frac 1{u_i-u_j}\\
\quad{}\times \det_{1\le i,j \le n}\! \left[(u_iq^{k_i})^{n-j}
\left(1-t_i^{j-n-1}\frac{1-t_0u_iq^{k_i}} {1-t_0u_iq^{k_i}/t_i}
\prod_{s=1}^n \frac{u_iq^{k_i}-u_s}
{u_iq^{k_i}/t_i-u_s}\right)\right]\;q^{\sum\limits_{i=1}^n(1-i)k_i}\\
\quad{}\times \prod_{i=1}^n \frac{(t_0u_iq/t_i,bu_i;q)_{k_i}}
{(t_0u_iq,t_0u_iq/dt_i;q)_{k_i}}t_i^{ik_i-\sum\limits_{j=1}^ik_j}
\cdot \frac{(d;q)_{|\mathbf k|}} {(t_0q/bt_1\cdots
t_n;q)_{|\mathbf k|}}
\left(\frac{t_0q}{bd}\right)^{|\mathbf k|}\\
\qquad{}=\frac{(t_0q/b,t_0q/bdt_1\cdots t_n;q)_\infty}
{(t_0q/bd,t_0q/bt_1\cdots t_n;q)_\infty}
\prod_{i=1}^n\frac{(t_0u_iq/t_i,t_0u_iq/d;q)_\infty}
{(t_0u_iq,t_0u_iq/dt_i;q)_\infty},
\end{gather*}
provided $|t_0q/bd|<1$.
\end{Corollary}

The following identity is obtained from Corollary~\ref{an65nt} by
letting $d\to q^{-M}$.

\begin{Corollary}[A multivariable
terminating ${}_6\phi_5$ summation]\label{an65n} Let $b$,
$t_0,t_1,\dots,t_n$, and $u_1,\dots,u_n$ be indeterminate, let $M$
be a nonnegative integer. Then
\begin{gather*}
\sum_{\underset{|{\mathbf k}|\le M}{k_1,\dots,k_n\ge 0}}
\prod_{i,j=1}^n\frac{{(qu_i/t_iu_j;q)}_{k_i}}
{{(qu_i/u_j;q)}_{k_i}}\, \prod_{1\le i<j\le
n}\frac{(t_ju_i/u_j;q)_{k_i-k_j}} {(qu_i/t_iu_j;q)_{k_i-k_j}}
\prod_{1\le i < j \le n}\frac 1{u_i-u_j}\\
\qquad{}\times \det_{1\le i,j \le n}\! \left[(u_iq^{k_i})^{n-j}
\left(1-t_i^{j-n-1}\frac{1-t_0u_iq^{k_i}} {1-t_0u_iq^{k_i}/t_i}
\prod_{s=1}^n \frac{u_iq^{k_i}-u_s}
{u_iq^{k_i}/t_i-u_s}\right)\right]\;q^{\sum\limits_{i=1}^n(1-i)k_i}\\
\qquad{}\times \prod_{i=1}^n \frac{(t_0u_iq/t_i,bu_i;q)_{k_i}}
{(t_0u_iq,t_0u_iq^{1+M}/t_i;q)_{k_i}}t_i^{ik_i-\sum\limits_{j=1}^ik_j}
\cdot\frac{(q^{-M};q)_{|\mathbf k|}} {(t_0q/bt_1\cdots
t_n;q)_{|\mathbf k|}} \left(\frac{t_0q^{1+M}}b\right)^{|\mathbf
k|}
\\
\qquad{}=\frac{(t_0q/b;q)_M} {(t_0q/bt_1\cdots t_n;q)_M}
\prod_{i=1}^n\frac{(t_0u_iq/t_i;q)_M} {(t_0u_iq;q)_M}.
\end{gather*}
\end{Corollary}

Another multivariable terminating ${}_6\phi_5$ summation theorem
is obtained from Corollary~\ref{an65nt} by letting $t_i\to
q^{m_i}$, $i=1,\dots,n$. We do not state it explicitly.

The last result in this section is, at this moment, only
conjectural (which we however verif\/ied numerically by computer).
We provide a heuristic (but not completely correct) derivation
using the $C_n$ ${}_8\phi_7$ summation in Proposition~\ref{cn87}
combined with the multivariable matrix inverse in
Corollary~\ref{cormi2c}.

\begin{Conjecture}[A multivariable
terminating ${}_8\phi_7$ summation]\label{cn87n} Let $a$, $b$,
$c$, $d$, and $u_1,\dots,u_n$ be indeterminate, let
$m_1,m_2,\dots,m_n$ be nonnegative integers. Then
\begin{gather}
\sum_{\underset{i=1,\dots,n}{0\le k_i\le m_i}}
\prod_{i,j=1}^n\frac{{(q^{1-m_i}u_i/u_j, au_iu_jq^{1+|\mathbf
m|-m_i};q)}_{k_i}\, (q^{m_j}u_i/u_j;q)_{k_i-k_j}}
{{(qu_i/u_j,au_iu_jq^{1+|\mathbf m|};q)}_{k_i}\,
(qu_i/u_j;q)_{k_i-k_j}}
\prod_{1\le i < j \le n}\frac{u_iq^{k_i}-u_jq^{k_j}}{u_i-u_j}\nonumber\\
\qquad{}\times \prod_{1\le i<j\le n}\frac{(au_iu_jq^{|\mathbf
m|};q)_{k_i+k_j}} {(au_iu_jq^{1+|\mathbf
m|-m_i-m_j};q)_{k_i+k_j}}\cdot
q^{|\mathbf k|+(n-1)\sum\limits_{i=1}^nm_ik_i}\nonumber\\
\qquad{}\times \prod_{1\le i < j \le
n}(u_i-u_j)^{-1}(1-q^{-|\mathbf m|}/au_iu_j)^{-1}
\prod_{i=1}^n(u_i+q^{-|\mathbf m|}/au_i)^{-1}\nonumber\\
\qquad{}\times \det_{1\le i,j \le n}\!
\bigg[(u_iq^{k_i}+q^{-k_i-m_i}/au_i)^{n+1-j}-
(u_iq^{k_i-m_i}+q^{-k_i+m_i-|\mathbf m|}/au_i)^{n+1-j}\nonumber\\
\qquad{}\times \prod_{s=1}^n
\frac{(u_iq^{k_i}-u_s)(1-q^{-k_i-|\mathbf m|}/au_iu_s)}
{(u_iq^{k_i-m_i}-u_s)(1-q^{-k_i+m_i-|\mathbf m|}/au_iu_s)}\bigg]\nonumber\\
\qquad{}\times \prod_{i=1}^n
\frac{(bu_i,cu_i,du_i,a^2u_iq^{1+|\mathbf m|}/bcd;q)_{k_i}}
{(au_iq^{1+|\mathbf m|-m_i}/b,au_iq^{1+|\mathbf m|-m_i}/c,
au_iq^{1+|\mathbf m|-m_i}/d,bcdu_iq^{-m_i}/a;q)_{k_i}}\nonumber\\
\qquad{}= \prod_{1\le i<j\le n}(au_iu_jq^{1+|\mathbf
m|-m_i-m_j};q)_{m_i+m_j}^{-1}
\prod_{i,j=1}^n(au_iu_jq^{1+|\mathbf m|-m_i};q)_{m_i}\nonumber\\
\qquad{}\times \frac{(aq/bc,aq/bd,aq/cd;q)_{|\mathbf m|}}
{\prod\limits_{i=1}^n(au_iq^{1+|\mathbf
m|-m_i}/b,au_iq^{1+|\mathbf m|-m_i}/c, au_iq^{1+|\mathbf
m|-m_i}/d,aq/bcdu_i;q)_{m_i}}.\label{cn87neq}
\end{gather}
\end{Conjecture}

{\em Heuristic derivation of Conjecture~{\rm \ref{cn87n}}.} We
apply multidimensional inverse relations. We have \eqref{rotinv2}
by the $u_i\mapsto u_iq^{l_i}$, $M_i\mapsto M_i-l_i$,
$i=1,\dots,n$, case of Proposition~\ref{cn87}, where
\begin{gather*}
b_{\mathbf k}=\!\prod_{1\le i < j \le n}\!
\frac{u_iq^{k_i}-u_jq^{k_j}}{u_i-u_j} \!\prod_{1\le i \le j \le
n}\! \frac{1-au_iu_jq^{k_i+k_j}}{1-au_iu_j} \!\prod_{i,j=1}^n
\!\frac{{(u_iq^{-M_j}/u_j;q)}_{k_i-k_j}\,(au_iu_j;q)_{k_i+k_j}}
{{(qu_i/u_j;q)}_{k_i-k_j}\,(au_iu_jq^{1+M_i};q)_{k_i+k_j}}\\
\phantom{b_{\mathbf k}=}{}\times
\prod_{i=1}^n\frac{(bu_i,cu_i,du_i,a^2u_iq^{1+|\mathbf
M|}/bcd;q)_{k_i}} {(au_iq/b,au_iq/c,au_iq/d,bcdu_iq^{-|\mathbf
M|}/a;q)_{k_i}} \cdot q^{|\mathbf k|}
\end{gather*}
and
\begin{gather*}
a_{\mathbf l}= \prod_{1\le i<j\le
n}(au_iu_jq^{1+l_i+l_j};q)_{M_i+M_j}^{-1}
\prod_{i,j=1}^n(au_iu_jq^{1+l_i+l_j};q)_{M_i}\\
\phantom{a_{\mathbf l}=}{}\times
\frac{(aq/bc,aq/bd,aq/cd;q)_{|\mathbf M|}}
{\prod\limits_{i=1}^n(au_iq^{1+l_i}/b,au_iq^{1+l_i}/c,
au_iq^{1+l_i}/d,aq^{1+|\mathbf M|-M_i-l_i}/bcdu_i;q)_{M_i}}\\
\phantom{a_{\mathbf l}=}{}\times \!\prod_{1\le i < j \le n}\!
\frac{u_iq^{l_i}-u_jq^{l_j}}{u_i-u_j} \!\prod_{1\le i \le j \le
n}\!\frac{1-au_iu_jq^{l_i+l_j}}{1-au_iu_j} \!\prod_{i,j=1}^n \!
\frac{{(u_iq^{-M_j}/u_j;q)}_{l_i-l_j}\,(au_iu_j;q)_{l_i+l_j}}
{{(qu_i/u_j;q)}_{l_i-l_j}\,(au_iu_jq^{1+M_i};q)_{l_i+l_j}}\\
\phantom{a_{\mathbf l}=}{}\times
\prod_{i=1}^n\frac{(bu_i,cu_i,du_i,a^2u_iq^{1+|\mathbf
M|}/bcd;q)_{l_i}} {(au_iq/b,au_iq/c,au_iq/d,bcdu_iq^{-|\mathbf
M|}/a;q)_{l_i}} \cdot q^{|\mathbf l|},
\end{gather*}
and $g_{\mathbf k\mathbf l}$ as in \eqref{cormigc} with
$t_i\mapsto q^{-M_i}$, $i=1,\dots,n$.

If we now consider \eqref{rotinv1} with the above sequences
$a_{\mathbf l}$ and $b_{\mathbf k}$, and $f_{\mathbf m\mathbf k}$
as in \eqref{cormifc} with $t_i\mapsto q^{-M_i}$, $i=1,\dots,n$,
but do everywhere the replacements $M_i\mapsto -M_i$,
$i=1,\dots,n$ (the latter would require justif\/ication), we
obtain, after the substitutions $a\mapsto aq^{|\mathbf m|}$,
$u_i\mapsto u_iq^{-k_i}$, $m_i\mapsto m_i+k_i$, $i=1,\dots,n$,
simplif\/ications and subsequent relabelling ($m_i\mapsto k_i$,
$M_i\mapsto m_i$, $i=1,\dots,n$), the identity in \eqref{cn87neq}.
\hfill\qed

We hope that Conjecture~\ref{cn87n} (assuming its correctness)
will be useful in  establishing explicit formulae related to
$B_n$, $C_n$ or $D_n$ Macdonald polynomials.

\section[Hypergeometric specialization of
$A_{n-1}$ Macdonald polynomials]{Hypergeometric specialization \\
of $\boldsymbol{A_{n-1}}$ Macdonald polynomials} \label{sechypsp}

We show here that the Pieri formula, Theorem~\ref{theopieri}, and
the (equivalent) recursion formula, Theorem~\ref{theomain}, both
can be viewed as multidimensional generalizations of the
terminating very-well-poised ${}_6\phi_5$ summation \eqref{sum65}.
This, together with the knowledge that the ${}_6\phi_5$ summation
has a nonterminating extension, see \eqref{sumnt65}, eventually
leads us in Section~\ref{secmacdrev} to extend the family of
Macdonald polynomials, indexed by partitions, to Macdonald
symmetric functions, indexed by $n$-tuples of complex numbers.

Following Macdonald~\cite[p.\ 338]{Ma}, we consider for an
indeterminate $u$ the homomorphism $\varepsilon_{u,t}:
\mathsf{Sym}\to\mathbb Q(q,t)$, def\/ined by
\begin{gather}\label{defevt}
\varepsilon_{u,t}(p_r)=\frac{1-u^r}{1-t^r},
\end{gather}
for each integer $r\ge 1$. In particular, if $u$ is replaced by
$t^n$ one has
\begin{gather*}
\varepsilon_{t^n,t}(p_r)=\frac{1-t^{nr}}{1-t^r}=
p_r(1,t,\dots,t^{n-1}),
\end{gather*}
and hence for any $f\in\mathsf{Sym}$
\begin{gather*}
\varepsilon_{t^n,t}(f)=f(1,t,\dots,t^{n-1}),
\end{gather*}
the ``principal specialization'' of $f$. (Another notation for the
principal specialization is $u_0(f)$~\cite[p.~331]{Ma}, while in
terms of plethystic notation it reads as $f[\frac{1-t^n}{1-t}]$;
more generally, $\varepsilon_{u,t}(f)=f[\frac{1-u}{1-t}]$,
relating to \eqref{defevt}.) Since (recall \eqref{mcdps})
\begin{gather}\label{prspecp}
\varepsilon_{t^n,t}P_{\la}(X;q,t)=t^{\sum\limits_{i=1}^n(i-1)\la_i}
\prod_{1\le i<j\le n}\frac{(t^{j-i+1};q)_{\la_i-\la_j}}
{(t^{j-i};q)_{\la_i-\la_j}},
\end{gather}
we have (see \cite[p.\ 338, Eq.~(6.17)]{Ma})
\begin{gather}\label{hypspecp}
\varepsilon_{u,t}P_{\la}(X;q,t)=t^{\sum\limits_{i=1}^n(i-1)\la_i}\,
\prod_{i=1}^n\frac{(ut^{1-i};q)_{\la_i}} {(t^{n+1-i};q)_{\la_i}}\,
\prod_{1\le i<j\le n}\frac{(t^{j-i+1};q)_{\la_i-\la_j}}
{(t^{j-i};q)_{\la_i-\la_j}},
\end{gather}
where $n\ge l(\la)$. (We can take any such $n$. The homomorphism
$\varepsilon_{u,t}$ in fact does not depend on~$n$.) This follows
from \eqref{prspecp} by a simple polynomial argument. Both sides
of \eqref{hypspecp} are polynomials in $u$ (with  coef\/f\/icients
in $\mathbb Q(q,t)$) which agree for inf\/initely many values of
$u$, namely whenever $u=t^n$ for any $n\ge l(\la)$, hence are
identically equal.

Combining \eqref{mcdnf}, \eqref{bqt} and \eqref{hypspecp}, we have
\begin{gather}
\varepsilon_{u,t}Q_{\la}(X;q,t)=
b_{\la}(q,t)\,\varepsilon_{u,t}P_{\la}(X;q,t)\nonumber\\
\phantom{\varepsilon_{u,t}Q_{\la}(X;q,t)}{}
=\prod_{i=1}^n\frac{(t^{n+1-i};q)_{\la_i}}{(qt^{n-i};q)_{\la_i}}\,
\prod_{1\le i<j\le n} \frac{(qt^{j-i},t^{j-i};q)_{\la_i-\la_j}}
{(t^{j-i+1},qt^{j-i-1};q)_{\la_i-\la_j}}\nonumber\\
\phantom{\varepsilon_{u,t}Q_{\la}(X;q,t)=}{}\times
t^{\sum\limits_{i=1}^n(i-1)\la_i}\,
\prod_{i=1}^n\frac{(ut^{1-i};q)_{\la_i}} {(t^{n+1-i};q)_{\la_i}}\,
\prod_{1\le i<j\le n}\frac{(t^{j-i+1};q)_{\la_i-\la_j}}
{(t^{j-i};q)_{\la_i-\la_j}}\nonumber\\
\phantom{\varepsilon_{u,t}Q_{\la}(X;q,t)}{}
=t^{\sum\limits_{i=1}^n(i-1)\la_i}\,
\prod_{i=1}^n\frac{(ut^{1-i};q)_{\la_i}} {(qt^{n-i};q)_{\la_i}}\,
\prod_{1\le i<j\le n}\frac{(qt^{j-i};q)_{\la_i-\la_j}}
{(qt^{j-i-1};q)_{\la_i-\la_j}},\label{hypspecq}
\end{gather}
where $n\ge l(\la)$.

In view of the explicit evaluations \eqref{hypspecp} and
\eqref{hypspecq}, it is clear how useful the application of the
homomorphism $\varepsilon_{u,t}$ is when applied to identities
involving Macdonald polynomials, in particular for deriving
multivariate formulae of hypergeometric type. In this respect
(when being applied to Macdonald polynomials), we may consider
$\varepsilon_{u,t}$ as a ``hypergeometric specialization'' (which
is slightly misleading as $\varepsilon_{u,t}$ is actually {\em
not} a specialization homomorphism for $u\neq t^n$; it can be
regarded as an analytic continued specialization homomorphism
though).

We mention that by applying $\varepsilon_{u,t}$ to the set $Y$ in
the Cauchy formula~\cite[p.\ 324, Eq.~(4.2)]{Ma},
\begin{gather*}
\sum_{\la}P_\la(X;q,t)Q_\la(Y;q,t)=
\prod_{i,j}\frac{(tx_iy_j;q)_\infty}{(x_iy_j;q)_\infty},
\end{gather*}
Macdonald~\cite[p.\ 374, Eq.~(4)]{Ma} obtained a $q$-binomial
theorem for multiple series of Macdonald polynomial argument. This
identity became in fact the starting point in the development of a
whole theory of identities for series of (classical $A_{n-1}$)
Macdonald polynomial argument, see Kaneko~\cite{Kan1,Kan2}, Baker
and Forrester~\cite{BF}, and Warnaar~\cite{W}.

\subsection{Hypergeometric specialization of the Pieri formula}
We apply $\varepsilon_{u,t}$ to both sides of the Pieri formula in
Theorem~\ref{theopieri}, utilizing \eqref{hypspecq}. On the
left-hand side we obtain
\begin{gather}
\varepsilon_{u,t}Q_{(\la_1,\dots,\la_n)}\,
\varepsilon_{u,t}Q_{(\la_{n+1})}=
t^{\sum\limits_{i=1}^n(i-1)\la_i}\,
\prod_{i=1}^n\frac{(ut^{1-i};q)_{\la_i}}
{(qt^{n-i};q)_{\la_i}}\nonumber\\
\phantom{\varepsilon_{u,t}Q_{(\la_1,\dots,\la_n)}\,
\varepsilon_{u,t}Q_{(\la_{n+1})}=}{}\times{} \prod_{1\le i<j\le
n}\frac{(qt^{j-i};q)_{\la_i-\la_j}}
{(qt^{j-i-1};q)_{\la_i-\la_j}}\cdot
\frac{(u;q)_{\la_{n+1}}}{(q;q)_{\la_{n+1}}}.\label{id1}
\end{gather}
Since
\begin{gather*}
\varepsilon_{u,t}Q_{(\la_1+\ta_1,\dots,\la_n+\ta_n,
\la_{n+1}-|\ta|)}=
t^{\big(\sum\limits_{i=1}^n(i-1)(\la_i+\ta_i)\big)+n(\la_{n+1}-|\ta|)}\\
\qquad{}\times \prod_{i=1}^n\frac{(ut^{1-i};q)_{\la_i+\ta_i}}
{(qt^{n+1-i};q)_{\la_i+\ta_i}}\cdot
\frac{(ut^{-n};q)_{\la_{n+1}-|\ta|}}
{(q;q)_{\la_{n+1}-|\ta|}}\\
\qquad{}\times \prod_{1\le i<j\le
n}\frac{(qt^{j-i};q)_{\la_i+\ta_i-\la_j-\ta_j}}
{(qt^{j-i-1};q)_{\la_i+\ta_i-\la_j-\ta_j}}\, \prod_{i=1}^n\frac
{(qt^{n+1-i};q)_{\la_i+\ta_i-\la_{n+1}+|\ta|}}
{(qt^{n-i};q)_{\la_i+\ta_i-\la_{n+1}+|\ta|}},
\end{gather*}
we obtain on the right-hand side
\begin{gather}
\sum_{\ta\in\mathbb N^n}
\prod_{i=1}^n\frac{(t,q^{|\ta|+1+\la_i-\la_{n+1}}t^{n-i};q)_{\ta_i}}
{(q,q^{|\ta|+\la_i-\la_{n+1}}t^{1+n-i};q)_{\ta_i}}\nonumber\\
\times \prod_{1\le i<j\le n}
\frac{(q^{\la_i-\la_j}t^{j-i+1},q^{-\ta_j+1+\la_i-\la_j}t^{j-i-1};q)_{\ta_i}}
{(q^{1+\la_i-\la_j}t^{j-i},q^{-\ta_j+\la_i-\la_j}t^{j-i};q)_{\ta_i}}\cdot
\varepsilon_{u,t}Q_{(\la_1+\ta_1,\dots,\la_n+\ta_n,
\la_{n+1}-|\ta|)}\nonumber\\
=t^{\sum\limits_{i=1}^{n+1}(i-1)\la_i}
\frac{(ut^{-n};q)_{\la_{n+1}}}{(q;q)_{\la_{n+1}}}
\prod_{i=1}^n\frac{(ut^{1-i};q)_{\la_i}}{(qt^{n+1-i};q)_{\la_i}}
\frac{(qt^{n+1-i};q)_{\la_i-\la_{n+1}}}{(qt^{n-i};q)_{\la_i-\la_{n+1}}}
\prod_{1\le i<j\le n}\frac{(qt^{j-i};q)_{\la_i-\la_j}}
{(qt^{j-i-1};q)_{\la_i-\la_j}}\nonumber\\
\times \sum_{\ta_1,\dots,\ta_n\ge 0} \prod_{1\le i<j\le
n}\frac{(1-q^{\la_i-\la_j+\ta_i-\ta_j}t^{j-i})}
{(1-q^{\la_i-\la_j}t^{j-i})}\,
\prod_{i=1}^n\frac{(1-q^{\la_i-\la_{n+1}+\ta_i+|\ta|}t^{n+1-i})}
{(1-q^{\la_i-\la_{n+1}}t^{n+1-i})}\nonumber\\
\times \prod_{i,j=1}^n \frac{(q^{\la_i-\la_j}t^{j-i+1};q)_{\ta_i}}
{(q^{1+\la_i-\la_j}t^{j-i};q)_{\ta_i}}\,
\prod_{i=1}^n\frac{(q^{\la_i-\la_{n+1}}t^{n+1-i};q)_{|\ta|}}
{(q^{1+\la_i-\la_{n+1}}t^{n-i};q)_{|\ta|}}\nonumber\\
\times \prod_{i=1}^n\frac{(q^{\la_i}ut^{1-i};q)_{\ta_i}}
{(q^{1+\la_i}t^{n+1-i};q)_{\ta_i}}\cdot
\frac{(q^{-\la_{n+1}};q)_{|\ta|}}{(q^{1-\la_{n+1}}t^n/u;q)_{|\ta|}}
\left(\frac
qu\right)^{|\ta|}q^{\sum\limits_{i=1}^n(i-1)\ta_i}.\label{id2}
\end{gather}
If we equate the right-hand side expressions in \eqref{id2} and
\eqref{id1}, divide both sides by
\begin{gather*}
t^{\sum\limits_{i=1}^{n+1}(i-1)\la_i}
\frac{(ut^{-n};q)_{\la_{n+1}}}{(q;q)_{\la_{n+1}}}
\prod_{i=1}^n\frac{(ut^{1-i};q)_{\la_i}}{(qt^{n+1-i};q)_{\la_i}}
\frac{(qt^{n+1-i};q)_{\la_i-\la_{n+1}}}{(qt^{n-i};q)_{\la_i-\la_{n+1}}}
\prod_{1\le i<j\le n}\frac{(qt^{j-i};q)_{\la_i-\la_j}}
{(qt^{j-i-1};q)_{\la_i-\la_j}}
\end{gather*}
and simplify, while rewriting the right-hand side using
\begin{gather*}
\frac{(u;q)_{\la_{n+1}}}{(ut^{-n};q)_{\la_{n+1}}}
\prod_{i=1}^n\frac{(q^{-\la_i}t^{i-n-1};q)_{\la_{n+1}}}
{(q^{-\la_i}t^{i-n};q)_{\la_{n+1}}}
=\frac{(q^{1-\la_{n+1}}/u;q)_{\la_{n+1}}}
{(q^{1-\la_{n+1}}t^n/u;q)_{\la_{n+1}}}
\prod_{i=1}^n\frac{(q^{1+\la_i-\la_{n+1}}t^{n+1-i};q)_{\la_{n+1}}}
{(q^{1+\la_i-\la_{n+1}}t^{n-i};q)_{\la_{n+1}}},
\end{gather*}
we obtain the following result:
\begin{gather*}
\sum_{\underset{|\ta|\le\la_{n+1}}{\ta_1,\dots,\ta_n\ge 0}}
\prod_{1\le i<j\le n}
\frac{(q^{\la_i+\ta_i}t^{-i}-q^{\la_j+\ta_j}t^{-j})}
{(q^{\la_i}t^{-i}-q^{\la_j}t^{-j})}\,
\prod_{i=1}^n\frac{(1-q^{\la_i-\la_{n+1}+\ta_i+|\ta|}t^{n+1-i})}
{(1-q^{\la_i-\la_{n+1}}t^{n+1-i})}\nonumber\\
\qquad{}\times \prod_{i,j=1}^n
\frac{(q^{\la_i-\la_j}t^{j-i+1};q)_{\ta_i}}
{(q^{1+\la_i-\la_j}t^{j-i};q)_{\ta_i}}\,
\prod_{i=1}^n\frac{(q^{\la_i-\la_{n+1}}t^{n+1-i};q)_{|\ta|}}
{(q^{1+\la_i-\la_{n+1}}t^{n-i};q)_{|\ta|}}\nonumber\\
\qquad{}\times \prod_{i=1}^n\frac{(q^{\la_i}ut^{1-i};q)_{\ta_i}}
{(q^{1+\la_i}t^{n+1-i};q)_{\ta_i}}\cdot
\frac{(q^{-\la_{n+1}};q)_{|\ta|}}{(q^{1-\la_{n+1}}t^n/u;q)_{|\ta|}}
\left(\frac qu\right)^{|\ta|}\nonumber\\
\qquad{}=\frac{(q^{1-\la_{n+1}}/u;q)_{\la_{n+1}}}
{(q^{1-\la_{n+1}}t^n/u;q)_{\la_{n+1}}}\,
\prod_{i=1}^n\frac{(q^{1+\la_i-\la_{n+1}}t^{n+1-i};q)_{\la_{n+1}}}
{(q^{1+\la_i-\la_{n+1}}t^{n-i};q)_{\la_{n+1}}}.
\end{gather*}
Now, this identity (which is equivalent to the classical
terminating ${}_6\phi_5$ summation when $n=1$) is simply the
special case of Milne's $A_{n-1}$ terminating ${}_6\phi_5$
summation in Proposition~\ref{an65}, where $a\mapsto t$, $b\mapsto
q^{\la_{n+1}}ut^{1-n}$, $c_i\mapsto t$, $u_i\mapsto
q^{\la_i-\la_{n+1}}t^{n-i}$ ($1\le i\le n$), and
$M\mapsto\la_{n+1}$. In particular, observe that the $n+1$
indeterminates $a,c_1,\dots,c_n$ all are substituted by $t$.

\subsection{Hypergeometric specialization of the recursion formula}

We apply $\varepsilon_{u,t}$ to both sides of the recursion
formula in Theorem~\ref{theomain}, utilizing \eqref{hypspecq}. On
the left-hand side we obtain
\begin{gather}\label{id4}
\varepsilon_{u,t}Q_{(\la_1,\dots,\la_{n+1})}=
t^{\sum\limits_{i=1}^{n+1}(i-1)\la_i}\,
\prod_{i=1}^{n+1}\frac{(ut^{1-i};q)_{\la_i}}
{(qt^{n+1-i};q)_{\la_i}}\, \prod_{1\le i<j\le
n+1}\frac{(qt^{j-i};q)_{\la_i-\la_j}}
{(qt^{j-i-1};q)_{\la_i-\la_j}}.
\end{gather}
Since
\begin{gather*}
\varepsilon_{u,t}Q_{(\la_{n+1}-|\ta|)}\,
\varepsilon_{u,t}Q_{(\la_1+\ta_1,\dots,\la_n+\ta_n)}=
\frac{(u;q)_{\la_{n+1}-|\ta|}}{(q;q)_{\la_{n+1}-|\ta|}}\cdot
t^{\sum\limits_{i=1}^n(i-1)(\la_i+\ta_i)}\\
\phantom{\varepsilon_{u,t}Q_{(\la_{n+1}-|\ta|)}\,
\varepsilon_{u,t}Q_{(\la_1+\ta_1,\dots,\la_n+\ta_n)}=}{}\times
\prod_{i=1}^n\frac{(ut^{1-i};q)_{\la_i+\ta_i}}
{(qt^{n-i};q)_{\la_i+\ta_i}}\, \prod_{1\le i<j\le
n}\frac{(qt^{j-i};q)_{\la_i+\ta_i-\la_j-\ta_j}}
{(qt^{j-i-1};q)_{\la_i+\ta_i-\la_j-\ta_j}},
\end{gather*}
we obtain on the right-hand side
\begin{gather}
\sum_{\ta\in\mathbb N^n} \prod_{i=1}^n t^{\ta_i}
\,\frac{(q/t,q^{1+\la_i-\la_{n+1}}t^{n-i};q)_{\ta_i}}
{(q,q^{1+\la_i-\la_{n+1}}t^{1+n-i};q)_{\ta_i}}\, \prod_{1\le i < j
\le n} \frac{{(q^{1+\la_i-\la_j}t^{j-i-1},
q^{-\ta_j+\la_i-\la_j}t^{j-i+1};q)}_{\ta_i}}
{{(q^{1+\la_i-\la_j}t^{j-i},q^{-\ta_j+\la_i-\la_j}t^{j-i};q)}_{\ta_i}}\nonumber\\
\qquad{}\times\prod_{1\le i < j \le n}\frac 1
{(q^{\la_i-\la_{n+1}+\ta_i}t^{n-i}-q^{\la_j-\la_{n+1}+\ta_j}t^{n-j})}
\det_{1\le i,j \le n}\!
\Bigg[(q^{\la_i-\la_{n+1}+\ta_i}t^{n-i})^{n-j}\nonumber\\
\left.\qquad{}\times
\left(1-t^{j-1}\frac{1-q^{\la_i-\la_{n+1}+\ta_i}t^{n-i+1}}
{1-q^{\la_i-\la_{n+1}+\ta_i}t^{n-i}} \prod_{s=1}^n
\frac{q^{\la_s-\la_{n+1}}t^{n-s}-
q^{\la_i-\la_{n+1}+\ta_i}t^{n-i}} {q^{\la_s-\la_{n+1}}t^{n-s+1}-
q^{\la_i-\la_{n+1}+\ta_i}t^{n-i}}\right)\right]\nonumber\\
\qquad{}\times \varepsilon_{u,t}Q_{(\la_{n+1}-|\ta|)}\,
\varepsilon_{u,t}Q_{(\la_1+\ta_1,\dots,\la_n+\ta_n)}\nonumber\\
\qquad{}
=t^{\sum\limits_{i=1}^n(i-1)\la_i}\,\frac{(u;q)_{\la_{n+1}}}{(q;q)_{\la_{n+1}}}
\prod_{i=1}^n\frac{(ut^{1-i};q)_{\la_i}} {(qt^{n-i};q)_{\la_i}}\,
\prod_{1\le i<j\le n}\frac{(qt^{j-i};q)_{\la_i-\la_j}}
{(qt^{j-i-1};q)_{\la_i-\la_j}}\nonumber\\\qquad{}\times
\sum_{\ta\in\mathbb N^n}
\prod_{i,j=1}^n\frac{{(q^{1+\la_i-\la_j}t^{j-i-1};q)}_{\ta_i}}
{{(q^{1+\la_i-\la_j}t^{j-i};q)}_{\ta_i}}\, \prod_{1\le i<j\le
n}\frac{(q^{\la_i-\la_j}t^{j-i+1};q)_{\ta_i-\ta_j}}
{(q^{1+\la_i-\la_j}t^{j-i-1};q)_{\ta_i-\ta_j}}\nonumber\\
\qquad{}\times \prod_{1\le i < j \le n}\frac 1
{(q^{\la_i-\la_{n+1}}t^{n-i}-q^{\la_j-\la_{n+1}}t^{n-j})}
\det_{1\le i,j \le n}\!
\Bigg[(q^{\la_i-\la_{n+1}+\ta_i}t^{n-i})^{n-j}\nonumber\\
\qquad{}\times
\left.\left(1-t^{j-1}\frac{1-q^{\la_i-\la_{n+1}+\ta_i}t^{n-i+1}}
{1-q^{\la_i-\la_{n+1}+\ta_i}t^{n-i}} \prod_{s=1}^n
\frac{q^{\la_s-\la_{n+1}}t^{n-s}-
q^{\la_i-\la_{n+1}+\ta_i}t^{n-i}} {q^{\la_s-\la_{n+1}}t^{n-s+1}-
q^{\la_i-\la_{n+1}+\ta_i}t^{n-i}}\right)\right]\label{id5}\\
\qquad{}\times\! \prod_{i=1}^n
\frac{(q^{1+\la_i-\la_{n+1}}t^{n-i},q^{\la_i}ut^{1-i};q)_{\ta_i}}
{(q^{1+\la_i-\la_{n+1}}t^{1+n-i},q^{1+\la_i}t^{n-i};q)_{\ta_i}}\cdot
\frac{(q^{-\la_{n+1}};q)_{|\ta|}}{(q^{1-\la_{n+1}}/u;q)_{|\ta|}}
\! \left(\frac {qt}u\right)^{|\ta|}\!\!
q^{\sum\limits_{i=1}^n(1-i)\ta_i}
t^{\sum\limits_{i=1}^n2(i-1)\ta_i}.\nonumber
\end{gather}
If we equate the right-hand side expressions in \eqref{id5} and
\eqref{id4}, divide both sides by
\begin{gather*}
t^{\sum\limits_{i=1}^n(i-1)\la_i}\,
\frac{(u;q)_{\la_{n+1}}}{(q;q)_{\la_{n+1}}}
\prod_{i=1}^n\frac{(ut^{1-i};q)_{\la_i}} {(qt^{n-i};q)_{\la_i}}\,
\prod_{1\le i<j\le n}\frac{(qt^{j-i};q)_{\la_i-\la_j}}
{(qt^{j-i-1};q)_{\la_i-\la_j}}
\end{gather*}
and simplify, while rewriting the right-hand side using
\begin{gather*}
t^{n\la_{n+1}}\frac{(ut^{-n};q)_{\la_{n+1}}}{(u;q)_{\la_{n+1}}}
\prod_{i=1}^n\frac{(qt^{n-i};q)_{\la_i}\,
(qt^{n+1-i};q)_{\la_i-\la_{n+1}}}
{(qt^{n+1-i};q)_{\la_i}\,(qt^{n-i};q)_{\la_i-\la_{n+1}}}\\
\qquad{}=\frac{(q^{1-\la_{n+1}}t^n/u;q)_{\la_{n+1}}}
{(q^{1-\la_{n+1}}/u;q)_{\la_{n+1}}}
\prod_{i=1}^n\frac{(q^{1+\la_i-\la_{n+1}}t^{n-i};q)_{\la_{n+1}}}
{(q^{1+\la_i-\la_{n+1}}t^{n+1-i};q)_{\la_{n+1}}},
\end{gather*}
we obtain the following result:
\begin{gather}
\sum_{\ta\in\mathbb N^n}
\prod_{i,j=1}^n\frac{{(q^{1+\la_i-\la_j}t^{j-i-1};q)}_{\ta_i}}
{{(q^{1+\la_i-\la_j}t^{j-i};q)}_{\ta_i}}\, \prod_{1\le i<j\le
n}\frac{(q^{\la_i-\la_j}t^{j-i+1};q)_{\ta_i-\ta_j}}
{(q^{1+\la_i-\la_j}t^{j-i-1};q)_{\ta_i-\ta_j}} \prod_{1\le i < j
\le n}\frac 1
{(q^{\la_i}t^{-i}-q^{\la_j}t^{-j})}\nonumber\\
\qquad{}\times \det_{1\le i,j \le n}\!
\left[(q^{\la_i+\ta_i}t^{-i})^{n-j}
\left(1-t^{j-1}\frac{1-q^{\la_i-\la_{n+1}+\ta_i}t^{n-i+1}}
{1-q^{\la_i-\la_{n+1}+\ta_i}t^{n-i}} \prod_{s=1}^n
\frac{q^{\la_s}t^{-s}- q^{\la_i+\ta_i}t^{-i}}
{q^{\la_s}t^{1-s}-q^{\la_i+\ta_i}t^{-i}}\right)\right]\nonumber\\
\qquad{}\times \prod_{i=1}^n
\frac{(q^{1+\la_i-\la_{n+1}}t^{n-i},q^{\la_i}ut^{1-i};q)_{\ta_i}}
{(q^{1+\la_i-\la_{n+1}}t^{1+n-i},q^{1+\la_i}t^{n-i};q)_{\ta_i}}\cdot
\frac{(q^{-\la_{n+1}};q)_{|\ta|}}{(q^{1-\la_{n+1}}/u;q)_{|\ta|}}
\!\left(\frac {qt}u\right)^{|\ta|}
\!q^{\sum\limits_{i=1}^n(1-i)\ta_i}
t^{\sum\limits_{i=1}^n2(i-1)\ta_i}\nonumber\\
\qquad{}= \frac{(q^{1-\la_{n+1}}t^n/u;q)_{\la_{n+1}}}
{(q^{1-\la_{n+1}}/u;q)_{\la_{n+1}}}
\prod_{i=1}^n\frac{(q^{1+\la_i-\la_{n+1}}t^{n-i};q)_{\la_{n+1}}}
{(q^{1+\la_i-\la_{n+1}}t^{n+1-i};q)_{\la_{n+1}}}.\label{id6}
\end{gather}
Now, this identity (which is equivalent to the classical
terminating ${}_6\phi_5$ summation when $n=1$) is simply the
special case of the multivariable terminating ${}_6\phi_5$
summation in Corollary~\ref{an65n}, where $a\mapsto t$, $b\mapsto
q^{\la_{n+1}}ut^{1-n}$, $c_i\mapsto t$, $u_i\mapsto
q^{\la_i-\la_{n+1}}t^{n-i}$ ($1\le i\le n$), and
$M\mapsto\la_{n+1}$. In particular, observe that the $n+1$
indeterminates $a,c_1,\dots,c_n$ all are substituted by $t$.

\section{More basic hypergeometric identities\\ involving
Macdonald polynomials} \label{secmorehyp}

In the previous section we showed that the Pieri formula and the
recursion formula, see Theo\-rems~\ref{theopieri} and
\ref{theomain}, both constitute (two dif\/ferent) multivariable
terminating ${}_6\phi_5$ summations which involve $A_{n-1}$
Macdonald polynomials. While nonterminating series are considered
in Section~\ref{secmacdrev}, one can ask whether other important
basic hypergeometric identities can be extended to the
multivariate setting involving $A_{n-1}$ Macdonald polynomials.
Concerning multivariate extensions of identities for
non-very-well-poised basic hypergeometric series, we refer to the
work of Kaneko~\cite{Kan1,Kan2}, Baker and Forrester~\cite{BF},
and Warnaar~\cite{W}, where several identities are established
that involve Macdonald polynomials playing the role of the {\em
argument} of the respective series. In the very-well-poised case,
which is investigated here, the $A_{n-1}$ Macdonald polynomials
play the role of {\em $q$-shifted factorials} (to which they would
reduce after principal specialization).

The main dif\/ference is the dimension of the series; in the
$A_{n-1}$ setting the very-well-poised structure is implicit from
taking the step $n\mapsto n+1$ to one higher dimension (see e.g.\
the proof of Proposition~\ref{an65}). Whereas for $n=1$ the
$A_{n-1}$ Macdonald polynomials are simply monomials, in the $n=2$
case they are, in view of
\begin{gather*}
P_{(\la_1,\la_2)}(x_1,x_2;q,t)
=(x_1x_2)^{\la_2}P_{(\la_1-\la_2,0)}(x_1,x_2;q,t)\\
\phantom{P_{(\la_1,\la_2)}(x_1,x_2;q,t)}{}
=(x_1x_2)^{\la_2}\frac{(q;q)_{\la_1-\la_2}}{(t;q)_{\la_1-\la_2}}
g_{\la_1-\la_2}(x_1,x_2;q,t)
\end{gather*}
(where we have used \eqref{fmacd} and \eqref{orc}), (multiples of
the) continuous $q$-ultraspherical polyno\-mials~$g_m$. The latter
are specif\/ic ${}_2\phi_1$ series, see \eqref{gm}, which when
principally specialized simplify to ratios of $q$-shifted
factorials by virtue of the Chu--Vandermonde summation theorem,
the terminating special case of the $q$-Gau{\ss} summation in
\eqref{sum21}.

We mention that for the nonreduced irreducible root system $BC_n$
very-well-poised basic hypergeometric series identities involving
Okounkov's~\cite{Ok} Macdonald interpolation polynomials or the
more general Koornwinder--Macdonald polynomials (both which are of
$BC_n$ type) have been established by Rains, see
\cite[Sec.~4]{Ra1} and \cite[Sec.~4]{Ra2}. (The $A_{n-1}$
identities investigated here appear to be essentially dif\/ferent
from related $BC_n$ identities found by Rains.)

One f\/irst question that arises is whether Jackson's terminating
very-well-poised ${}_8\phi_7$ summation \eqref{sum87} can be
extended to the multivariate setting involving $A_{n-1}$ Macdonald
polynomials. The answer is af\/f\/irmative. Consider the Pieri
formula in Theorem~\ref{theopieri}. It reads as
\begin{gather}\label{p1}
Q_{(\la_1,\ldots,\la_n)} \: Q_{(\la_{n+1})}= \sum_{\ta\in
\mathbb{N}^n} d_{\ta_1,\ldots,\ta_n}^{(q,t)} (u_1,\ldots,u_n)\:
Q_{(\la_1+\ta_1,\ldots,\la_n+\ta_n,\la_{n+1}-|\ta|)},
\end{gather}
with coef\/f\/icients $d_{\ta_1,\ldots,\ta_n} (u_1,\ldots,u_n)$
def\/ined in Subsection~\ref{secpieri}. We already know that
equation~\eqref{p1}, as it stands, represents a multivariable
terminating very-well-poised ${}_6\phi_5$ summation. The ``trick''
now is to restrict the set $X$ to $n$ variables, i.e.\ to take
$|X|=n$. Due to property~\eqref{stab} the Macdonald polynomials
appearing in the sum all vanish unless $|\ta|=\la_{n+1}$. We thus
obtain
\begin{gather}\label{p2}
Q_{(\la_1,\ldots,\la_n)} \: Q_{(\la_{n+1})}= \sum_{\ta\in
\mathbb{N}^n,|\ta|=\la_{n+1}} d_{\ta_1,\ldots,\ta_n}^{(q,t)}
(u_1,\ldots,u_n)\: Q_{(\la_1+\ta_1,\ldots,\la_n+\ta_n)},
\end{gather}
under the assumption $X=\{x_1,\dots,x_n\}$. It turns out that when
$\varepsilon_{u,t}$ is applied to both sides of the identity in
\eqref{p2} (in Section~\ref{sechypsp} we coined this as
``hypergeometric specialization''), its $n=2$ case reduces to a
variant of Jackson's ${}_8\phi_7$ summation \eqref{sum87}. In
fact, if we replace $n$ by $n+1$ in \eqref{p2}, and substitute
$\ta_{n+1}$ by $\la_{n+2}-(\ta_1+\cdots+\ta_n)$, we obtain an
identity which, after application of $\varepsilon_{u,t}$, is
essentially equivalent to Milne's $A_{n-1}$ extension of Jackson's
${}_8\phi_7$ summation in Proposition~\ref{an87}. Therefore, we
regard \eqref{p2} as a multivariable extension of Jackson's
very-well-poised ${}_8\phi_7$ summation.

On the other hand, we were (so far) not able to deduce a
multivariable ${}_8\phi_7$ summation directly from the
(nonspecialized) recursion formula in Theorem~\ref{theomain}.

It would be interesting to f\/ind a multivariable extension of
Bailey's terminating very-well-poised ${}_{10}\phi_9$
transformation formula (cf.\ \cite[Eq.~(III.28)]{GR}) involving
$A_{n-1}$ Macdonald polyno\-mials.

\section{Macdonald symmetric functions indexed by partitions\\ with
complex parts} \label{secmacdrev}

In this section we shall assume $|q|<1$. We use the recursion in
Theorem~\ref{theomain} now to {\em define} Macdonald symmetric
functions $Q_\la$ when $\la=(\la_1,\dots,\la_n)$ is {\em any}
sequence of complex numbers. One dif\/f\/iculty is to properly
def\/ine the one row case. (Another issue is convergence, since we
will be considering nonterminating sums.)

Kadell~\cite{Ka} used the classical def\/inition of a Schur
function in terms of a ratio of alternants to extend Schur
functions to partitions with complex parts. (Independently,
Danilov and Koshevoy~\cite{DK} def\/ine ``continuous Schur
functions'' by a multidimensional integral, with respect to 
a~Lebesgue measure in $\mathbb R^{n(n-1)/2}$, over all points of a
particular polytope, and show by an inductive argument that these
functions generalize the ratio of alternants formula for Schur
functions. In fact, Danilov and Koshevoy's continuous Schur
functions correspond exactly to Kadell's Schur functions indexed
by partitions with {\em real} parts.) We want to stress that our
proposed extension of Macdonald polynomials to complex parts (see
below) does {\em not} reduce to Kadell's extension of Schur
functions when $q=t$.

We shall begin with a f\/inite number of variables, say
$X=\{x_0,x_1,\dots,x_r\}$. (For convenience, we start to label $X$
with $0$ here). First, consider $m$ to be a nonnegative integer.
By appealing to the $q$-binomial theorem in \eqref{sum10} it
follows from taking coef\/f\/icients of $u^m$ in the generating
function in \eqref{genserg} that the one row Macdonald polynomials
$Q_{(m)}(X;q,t)=g_m(X;q,t)$ can be written in the following
explicit form:
\begin{gather*}
g_m=\underset{k_0+\dots+k_r=m}{\sum_{k_0,\dots,k_r\ge 0}}
\prod_{i=0}^r\frac{(t;q)_{k_i}}{(q;q)_{k_i}}\,x_i^{k_i}\\
\phantom{g_m}{}= \underset{0\le k_1+\dots+k_r\le
m}{\sum_{k_1,\dots,k_r\ge 0}}
\frac{(t;q)_{m-(k_1+\dots+k_r)}}{(q;q)_{m-(k_1+\dots+k_r)}}\,
x_0^{m-(k_1+\dots+k_r)}\,
\prod_{i=1}^r\frac{(t;q)_{k_i}}{(q;q)_{k_i}}x_i^{k_i}.
\end{gather*}
Although we do not need it here, we mention that for $r=1$ and
$x_0x_1=1$, the $g_m$ reduce to the continuous $q$-ultraspherical
polynomials of degree $m$ in the argument $(x_0+x_1)/2$,
considered in \cite[Ex.~1.29]{GR}.

We rewrite the above expression yet further, using the short
notation $|{\mathbf k}|=k_1+\dots+k_r$, and obtain the following
explicit form:
\begin{gather}\label{gm}
g_m=\frac{(t;q)_m}{(q;q)_m}\,x_0^m \underset{0\le |{\mathbf k}|\le
m}{\sum_{k_1,\dots,k_r\ge 0}} \frac{(q^{-m};q)_{|{\mathbf
k}|}}{(q^{1-m}/t;q)_{|{\mathbf k}|}}
\prod_{i=1}^r\frac{(t;q)_{k_i}}{(q;q)_{k_i}}
\left(\frac{qx_i}{tx_0}\right)^{k_i}.
\end{gather}

Using the def\/inition
\begin{gather*}
(a;q)_c=\frac{(a;q)_\infty}{(aq^c;q)_\infty}
\end{gather*}
(recall $|q|<1$) for any complex number $c$, we propose the
following def\/inition for a one row complex Macdonald
function\footnote{We call these new objects {\em functions} as
they are not anymore polynomials. However, note that the term {\em
Macdonald functions} is also used to denote other objects, namely
{\em modified Bessel functions of the second kind}~\cite{Wt}. (The
latter have nothing to do with algebraist I.G.~Macdonald who
introduced the $P_\la(X;q,t)$.)}:{\samepage
\begin{gather}\label{gc}
Q_{(c)}=g_c=\frac{(tx_0;q)_c}{(q;q)_c}
\frac{(q/tx_0;q)_{-c}}{(q/t;q)_{-c}} \sum_{k_1,\dots,k_r\ge 0}
\frac{(q^{-c};q)_{|{\mathbf k}|}}{(q^{1-c}/t;q)_{|{\mathbf k}|}}
\prod_{i=1}^r\frac{(t;q)_{k_i}}{(q;q)_{k_i}}
\left(\frac{qx_i}{tx_0}\right)^{k_i},
\end{gather}
which converges (if the series does not terminate) for
$|qx_i/tx_0|<1$ ($1\le i\le r$).}

We emphasize that \eqref{gc} is {\em not} an analytic continuation
of \eqref{gm}. In fact, \eqref{gc} is neither analytic in $c$, nor
in $q^c$ (in any domain around the origin). Another problem is
that $g_{c}$ is not symmetric in all the $x_i$ ($0\le i\le r$) but
only in the last $r$ of the $x_i$ ($1\le i\le r$). Indeed, already
for $r=1$, the $_2\phi_1$ transformation
$g_c(x_0,x_1)=g_c(x_1,x_0)$ turns out to be false if $c$ is not an
integer.

There are other possibilities to extend $g_m$ to complex numbers.
By our def\/inition \eqref{gc}, if $c$ is not a nonnegative
integer and $q=t$ (the Schur function case), then we get
$g_c=\infty$, which is dif\/ferent from
Kadell's~\cite[Eq.~(2.1)]{Ka} $s_{(c)}(x_0)=x_0^c=e^{c\ln(x_0)}$.
On the other hand, if we would have def\/ined $g_c$ by \eqref{gm}
(for complex $m=c$; relaxing the restriction $|{\mathbf k}|\le m$
of summation), our def\/inition would have also not matched
Kadell's in the $q=t$ case since after letting $q=t$ we would be
left with a product of geometric series on the right-hand side.
However, our particular choice of \eqref{gc} is motivated by some
nice properties, among which are \eqref{psr} and \eqref{qx}.

Since $g_c(x_0,x_1,\dots,x_r,0)=g_c(x_0,x_1,\dots,x_r)$, we may
let $r\to\infty$ (compare to \cite[p.~41]{Ma}). In the following,
we relax the restriction of $X$ being f\/inite. Thus, we allow
$r\in\mathbb N\cup\infty$.

Having provided a def\/inition of $Q_{(c)}$ for any complex number
$c$, it is now straightforward to extend Theorem~\ref{theomain} to
Macdonald functions indexed by sequences of complex numbers. Let
$\la=(\la_1,\dots ,\la_{n+1})$ be an arbitrary sequence of complex
numbers. We do {\em not} require $n+1\le |X|$. For any $1\le i \le
n+1$ def\/ine $u_i=q^{\la_i-\la_{n+1}}t^{n-i}$. Then
$Q_\la(X;q,t)$ is def\/ined recursively by \eqref{gc} and
\begin{gather}\label{qcn}
Q_{(\la_1,\ldots,\la_{n+1})}= \sum_{\ta\in\mathbb{N}^n}
c^{(q,t)}_{\ta_1,\ldots,\ta_n} (u_1,\ldots,u_n)\:
Q_{(\la_{n+1}-|\ta|)} \: Q_{(\la_1+\ta_1,\ldots,\la_n+\ta_n)},
\end{gather}
where $c^{(q,t)}_{\ta_1,\ldots,\ta_n} (u_1,\ldots,u_n)$ is the
same as in Subsection~\ref{secmain}.

While the (f\/inite) recursion in Theorem~\ref{theomain} was
proved by inverting the known Pieri formula for Macdonald
polynomials indexed by partitions, \eqref{qcn} now {\em defines}
Macdonald functions in the general case. The expansion in
\eqref{qcn} is in general inf\/inite and converges (when it does
not terminate) for $|q|, |qx_i/tx_0|<1$ ($i\ge 1$). As a matter of
formal manipulations (using multidimensional inverse relations),
the equivalence of Theorems~\ref{theomain} and \ref{theopieri}
(Pieri formula) for these Macdonald functions of complex parts is
immediate.

We do not know whether the complex $Q_\la$ form a family of
orthogonal functions, nor whether they are eigenfunctions of the
Macdonald operator (or some reasonable extension of this
operator). These questions, among others, wait for investigation.
What at all makes these extended objects interesting, then? In
fact, it turns out that for these ``complex Macdonald functions''
a~generalization of (at least) one of the so-called Macdonald
(ex-)conjectures holds; in particular, they satisfy an explicit
{\em evaluation formula}.

Observe that we have departed from the algebra of symmetric
functions in $X=\{x_0,x_1,\dots\}$. We are working in another
algebra. It is necessary to provide some details.

It is convenient to make the following def\/initions. For a
complex number $c$, introduce the following complex
``$q,t$-powers'' of $x$:
\begin{gather}\label{qtpowers}
x^{[c]}=x^{[c;q,t]}=\frac{(tx;q)_c}{(t;q)_c}
\frac{(q/tx;q)_{-c}}{(q/t;q)_{-c}}.
\end{gather}
Note that if $k$ is an integer then $x^{[k]}=x^k$. More generally,
$x^{[c+k]}=x^{[c]}x^k$.

Next, using ``$q,t$-powers'', extend the def\/inition of power
sums to complex numbers $c$:
\begin{gather*}
p_c(x_0,x_1,\dots)=\sum_{i\ge 0}x_i^{[c]}.
\end{gather*}
(We will actually only need the one-variable case
$p_c(x_0)=x_0^{[c]}$ here.) As usual, this def\/inition may be
extended to multiindices, $p_{(c_1,\dots,c_n)}=p_{c_1}\cdots
p_{c_n}$.

Let $\mathbb C$ denote the set of complex numbers. The algebra we
are considering is (algebraically) generated by the uncountably
inf\/inite set of products
\begin{gather*}
\{p_c(x_0)\,p_s(x_1,x_2,\dots)|\: c\in\mathbb C, s\in\mathbb
N-{0}\}
\end{gather*}
with coef\/f\/icients in $\mathbb C((q,t))$. Note that the above
set of products is not an algebraic basis as we do not have
algebraic independence (e.g.,
$(p_{c_1}p_{s_1})(p_{c_2}p_{s_2})=(p_{c_1}p_{s_2})(p_{c_2}p_{s_1})$).
We denote this algebra by $\widehat{\mathsf{Sym}}$.

There is a $\widehat{\mathsf{Sym}}$-extension of the homomorphism
$\varepsilon_{u,t}$, def\/ined in \cite[p.~338, Eq.~(6.16)]{Ma},
which acts nicely on $Q_{\la}$. For an indeterminate $u$, def\/ine
a homomorphism $\widehat\varepsilon_{u,t}:
\widehat{\mathsf{Sym}}\to C((q,t))$ by
\begin{gather*}
\widehat\varepsilon_{ut,t}\big[p_c(x_0)p_s(x_1,x_2,\dots)\big]=
u^{[c]}\,\frac{1-u^s}{1-t^s}
\end{gather*}
for each complex number $c$ and positive integer $s$. In
particular, if $u$ is replaced by $t^r$, we have
\begin{gather*}
\widehat\varepsilon_{t^{r+1},t}\big[p_c(x_0)p_s(x_1,x_2,\dots)\big]=
(t^r)^{[c]}\,\frac{1-t^{rs}}{1-t^s}=p_c(t^r)p_s(t^{r-1},\dots,t,1)
\end{gather*}
and hence for any $f\in\widehat{\mathsf{Sym}}$
\begin{gather*}
\widehat\varepsilon_{t^{r+1},t}(f)=f(t^r,t^{r-1},\dots,t,1).
\end{gather*}
(Compare this to the usual
$\varepsilon_{t^{r+1},t}(f)=f(t^r,t^{r-1},\dots,t,1)$, for any
$f\in\mathsf{Sym}$.)

We now have the {\em evaluation formula}
\begin{gather}\label{psr}
\widehat\varepsilon_{u,t}Q_{(\la_1,\dots,\la_n)}=
\prod_{i=1}^n\frac{(u;q)_{\la_i}}{(qt^{n-i};q)_{\la_i}}
\frac{(q/u;q)_{-\la_i}}{(qt^{i-1}/u;q)_{-\la_i}} \prod_{1\le
i<j\le n}\frac{(qt^{j-i};q)_{\la_i-\la_j}}
{(qt^{j-i-1};q)_{\la_i-\la_j}},
\end{gather}
where $\la_i\in\mathbb C$. We will present a pure basic
hypergeometric proof (which reduces to a new proof of the usual
analytic continued principal specialization formula for $Q_\la$ if
$\la$ is a partition).

We proceed by induction on $n$. For $n=1$ we f\/irst consider
$\widehat\varepsilon_{t^{r+1},t}Q_{(c)}$ with $Q_{(c)}$ given in
\eqref{gc}. By the $a_i\mapsto t$, $u_i\mapsto u^i$, $i=1,\dots,
n$, $b\mapsto q^{-c}$, and $c\mapsto q^{1-c}/t$ case of
Proposition~\ref{annt21}, we see that
\begin{gather*}
\widehat\varepsilon_{t^{r+1},t}Q_{(c)}=\frac{(t^{r+1};q)_c}{(q;q)_c}
\frac{(q/t^{r+1};q)_{-c}}{(q/t;q)_{-c}} \sum_{k_1,\dots,k_r\ge 0}
\frac{(q^{-c};q)_{|{\mathbf k}|}}{(q^{1-c}/t;q)_{|{\mathbf k}|}}
\prod_{i=1}^r\frac{(t;q)_{k_i}}{(q;q)_{k_i}}
\left(qt^{-1-i}\right)^{k_i}\\
\phantom{\widehat\varepsilon_{t^{r+1},t}Q_{(c)}}{}=
\frac{(t^{r+1};q)_c}{(q;q)_c}
\frac{(q/t^{r+1};q)_{-c}}{(q/t;q)_{-c}}\,
\frac{(q^{1-c}/t^{r+1},q/t;q)_\infty}{(q^{1-c}/t,q/t^{r+1};q)_\infty}=
\frac{(t^{r+1};q)_c}{(q;q)_c},
\end{gather*}
which is exactly the $u=t^{r+1}$ case of the right-hand side of
\eqref{psr} for $n=1$. Since this holds for $r=0,1,2,\dots$, by
analytic continuation we may replace $t^{r+1}$ by $u$, which
establishes the $n=1$ case of \eqref{psr}. For the inductive step,
we assume \eqref{psr} for partitions $\la$ with $l(\la)\le n$.
Apply~$\widehat\varepsilon_{u,t}$ to both sides of the recursion
formula \eqref{qcn} and use the inductive hypothesis to simplify
the summand. We are done if we can show that the sum evaluates to
\begin{gather*}
\prod_{i=1}^{n+1}\frac{(u;q)_{\la_i}}{(qt^{n+1-i};q)_{\la_i}}
\frac{(q/u;q)_{-\la_i}}{(qt^{i-1}/u;q)_{-\la_i}} \prod_{1\le
i<j\le n+1}\frac{(qt^{j-i};q)_{\la_i-\la_j}}
{(qt^{j-i-1};q)_{\la_i-\la_j}}.
\end{gather*}
This, however, follows by an application of the multivariable
${}_6\phi_5$ summation in Corollary~\ref{an65nt} (after performing
the substitutions $a\mapsto t$, $b\mapsto q^{\la_{n+1}}ut^{1-n}$,
$c_i\mapsto t$, $u_i\mapsto q^{\la_i-\la_{n+1}}t^{n-i}$, $1\le
i\le n$, and $M\mapsto\la_{n+1}$; compare with
equation~\eqref{id6}).\hfill\qed

\smallskip

There is a well-known duality formula for Macdonald polynomials
(cf.\ \cite[p.~332, Eq.~(6.6)]{Ma}),
\begin{gather*}
\frac{Q_\la(q^{\mu}t^{\delta})}{Q_\la(t^{\delta})}=
\frac{Q_\mu(q^{\la}t^{\delta})}{Q_\mu(t^{\delta})},
\end{gather*}
for partitions $\la$ and $\mu$ of length $\le n$, where
\begin{gather*}
Q_\la(q^{\mu}t^{\delta})=
Q_\la(q^{\mu_1}t^{n-1},q^{\mu_2}t^{n-2},\dots,q^{\mu_n}).
\end{gather*}

We do not know (at present) whether this relation still holds for
arbitrary complex sequences $\la_1,\la_2,\dots$, and
$\mu_1,\mu_2,\dots$. However, it does hold if the length of the
partitions is one. Namely, if $\la=(c)$, $\mu=(d)$ are one row
complex partitions, we have
\begin{gather*}
\frac{Q_{(c)}(q^dt,1)}{Q_{(c)}(t,1)}= \frac{(q;q)_c}{(t^2;q)_c}\,
\frac{(tx_0;q)_c}{(q;q)_c} \frac{(q/tx_0;q)_{-c}}{(q/t;q)_{-c}}
\sum_{k\ge 0} \frac{(q^{-c},t;q)_k}{(q^{1-c}/t,q;q)_k}
\left(\frac{q^{1-d}}{t^2}\right)^k\\
\phantom{\frac{Q_{(c)}(q^dt,1)}{Q_{(c)}(t,1)}}{}
=\frac{(tx_0;q)_c}{(t^2;q)_c}
\frac{(q/tx_0;q)_{-c}}{(q/t;q)_{-c}}\
{}_2\phi_1\!\left[\begin{matrix}q^{-c},t\\
q^{1-c}/t\end{matrix}\,;q,q^{1-d}/t^2\right].
\end{gather*}
The duality is now an immediate consequence of the iterate of
Heine's transformation \cite[Eq.~(III.2)]{GR},
\begin{gather*}
{}_2\phi_1\!\left[\begin{matrix}a,b\\
c\end{matrix}\,;q,z\right]=\frac{(c/b,bz;q)_\infty}{(c,z;q)_\infty}\;
{}_2\phi_1\!\left[\begin{matrix}abz/c,b\\
bz\end{matrix}\,;q,\frac cb\right],
\end{gather*}
valid for $\max(|z|,|c/b|)<1$.

To prevent possible misconception, we note that the well-known
property valid for Macdonald polynomials indexed by partitions,
$Q_\la(x_1,\dots,x_r)=0$ if $l(\la)>r$ (see \eqref{mcdnf} and
\eqref{stab}), does {\em not} hold in the general complex case.
For instance, if $X=\{x\}$ contains only one variable, then
\begin{gather}\label{qx}
Q_{(\la_1,\dots,\la_n)}(x)=
\prod_{i=1}^n\frac{(tx;q)_{\la_i}}{(qt^{n-i};q)_{\la_i}}
\frac{(q/tx;q)_{-\la_i}}{(qt^{i-2};q)_{-\la_i}} \prod_{1\le i<j\le
n}\frac{(qt^{j-i};q)_{\la_i-\la_j}}
{(qt^{j-i-1};q)_{\la_i-\la_j}},
\end{gather}
where $\la_i\in\mathbb C$. This formula (which can be proved by
induction, similar to the above proof of~\eqref{psr}) is indeed
independent from the representation of $\la$, i.e.\ we may choose
$(\la_1,\dots,\la_n)=(\la_1,\dots,\la_n,0,\dots,0)$, adding an
arbitrary number of zeros at the end of sequence. It is clear from
\eqref{qx} that if $\la$ is a usual integer partition, then
\begin{gather*}
Q_{(\la_1,\dots,\la_n)}(x)=\frac{(t;q)_{\la_1}}{(q;q)_{\la_1}}\,
x^{\la_1}\,\delta_{\la_20}\dots\delta_{\la_n0}
\end{gather*}
(where we were using \eqref{qtpowers} and $x^{[k]}=x^k$ for
integer $k$).

In this section, we extended the Macdonald polynomials $Q_\la$ to
arbitrary sequences $\la=(\la_1,\dots,\la_n)$ of complex numbers.
To give such an extension for $P_\la$ one may simply invoke
$Q_\la=b_\la P_\la$, see \eqref{mcdnf}, with the known explicit
expression of $b_\la=b_\la(q,t)$, extended to complex sequences
$\la$. (This does not mean that we necessarily assume
$b_\la=\langle Q_\la,Q_\la\rangle_{q,t}$ beforehand. Nevertheless,
the latter equality should conjecturally still hold, for some
suitable $\widehat{\mathsf{Sym}}$-extension of the inner product
$\protect{\langle\,,\,\rangle_{q,t}}$.) To utilize
Theorem~\ref{theodual} to def\/ine $P_\la$ in the complex case
does not make sense since the indexing partitions are given there
in the form $\la=(1^{m_1},2^{m_2},\dots$, $(n+1)^{m_{n+1}})$, and the
multiplicities $m_i$ have no meaning for partitions with complex
parts. (Already for compositions $\la$ one would run into trouble
here.)

Again, we do not yet know whether the Macdonald functions for
partitions with complex parts form a family of orthogonal
functions. They very well may be orthogonal, possibly with respect
to a $\widehat{\mathsf{Sym}}$ variant of the inner product in
\cite[p.~372, Eq.~(9.10)]{Ma}. In view of \eqref{psr}
and~\eqref{qx}, some nice properties do exist, which provides some
evidence that these new objects merit further investigation.

\begin{appendix}

\section{A multidimensional matrix inverse}\label{secmmi}

Here we state the multidimensional matrix inverse of
\cite[Th.~2.6]{LS}, which happens to be so crucial for the results
in this paper, and work out some special cases we need which we
give as corollaries.

Let $\mathbb Z$ be the set of integers and $n$ some positive
integer. For multi-integers $\mathbf m,\mathbf k\in\mathbb Z^n$,
we write $\mathbf m\ge\mathbf k$ for $m_i\ge k_i$, for
$i=1,\dots,n$.

We say that an inf\/inite $n$-dimensional matrix $F=(f_{\mathbf
m\mathbf k})_{\mathbf m,\mathbf k\in\mathbb Z^n}$ is
lower-triangular if $f_{\mathbf m\mathbf k}=0$ unless $\mathbf
m\ge\mathbf k$. When all $f_{\mathbf k\mathbf k}\neq 0$, there
exists a unique lower-triangular matrix $G=(g_{\mathbf m\mathbf
k})_{\mathbf m,\mathbf k\in\mathbb Z^n}$, called the inverse of
$F$, such that the following orthogonality relation holds:
\begin{gather}\label{orthrel}
\sum_{\mathbf m\ge\mathbf k}f_{\mathbf m\mathbf k}g_{\mathbf
k\mathbf l} =\delta_{\mathbf m\mathbf l},\qquad \text{for
all}\quad \mathbf m,\mathbf l\in\mathbb Z^n.
\end{gather}
Clearly, as $F$ and $G$ are both lower-triangular, $F$ must then
also be the inverse of $G$ and the dual relation,
\begin{gather}\label{orthreld}
\sum_{\mathbf m\ge\mathbf k}g_{\mathbf m\mathbf k}f_{\mathbf
k\mathbf l} =\delta_{\mathbf m\mathbf l},\qquad \text{for
all}\quad \mathbf m,\mathbf l\in\mathbb Z^n,
\end{gather}
must hold at the same time. Therefore, if $F$ and $G$ are
inf\/inite lower-triangular $n$-dimensional matrices and one of
the relations \eqref{orthrel} or \eqref{orthreld} hold, we say
that $F$ and $G$ are inverses of each other, or simply that $F$
and $G$ are mutually inverse.

It is immediate from the orthogonality relations \eqref{orthrel}
and \eqref{orthreld} that if $F$ and $G$ are mutually inverse, the
following two equations (a.k.a.\ ``inverse relations''):
\begin{subequations}\label{invrel}
\begin{gather}\label{invrel1}
\sum_{\mathbf k\in\mathbb Z^n}f_{\mathbf m\mathbf k}a_{\mathbf k}
=b_{\mathbf m},\qquad \text{for all} \quad {\mathbf m},
\end{gather}
and
\begin{gather}\label{invrel2}
\quad \sum_{\mathbf l\in\mathbb Z^n}g_{\mathbf k\mathbf
l}b_{\mathbf l} =a_{\mathbf k},\qquad \text{for all}\quad {\mathbf
k}
\end{gather}
\end{subequations}
are equivalent, provided both sides terminate or converge.

Similarly, the following two equations, where one sums over the
{\em first} (instead of the second) multi-index of the respective
matrices, are equivalent:
\begin{subequations}\label{rotinv}
\begin{gather}\label{rotinv1}
\sum_{\mathbf m\in\mathbb Z^n}f_{\mathbf m\mathbf k}a_{\mathbf m}
=b_{\mathbf k},\qquad \text{for all}\quad {\mathbf k},
\end{gather}
and
\begin{gather}\label{rotinv2}
\quad \sum_{\mathbf k\in\mathbb Z^n}g_{\mathbf k\mathbf
l}b_{\mathbf k} =a_{\mathbf l},\qquad \text{for all}\quad \mathbf
l,
\end{gather}
\end{subequations}
again, provided both sides terminate or converge.

The inverse relations \eqref{invrel} or \eqref{rotinv} are useful
tools for proving identities. For instance, if~\eqref{rotinv2}
holds for the inf\/inite lower-triangular matrix $(g_{\mathbf
m\mathbf k})_{\mathbf m,\mathbf k\in\mathbb Z^n}$ and inf\/inite
sequences $(b_{\mathbf k})_{\mathbf k\in\mathbb Z^n}$ and
$(a_{\mathbf l})_{\mathbf l\in\mathbb Z^n}$, then (assuming one
knows $F$, the inverse of $G$), the equation \eqref{rotinv1} must
automatically hold (subject to convergence). It is exactly this
form of inverse relations which is being utilized in this paper.

The following matrix inversion was derived in \cite[Th.~2.6]{LS}.

\begin{Proposition}\label{thmi2}
Let $b$ be an indeterminate and $a_i(k)$, $c_i(k)$ $(k\in \mathbb
Z,\ 1\le i \le n)$ be arbitrary sequences of indeterminates. Then
the infinite lower-triangular $n$-dimensional matrices
$(f_{\mathbf{m k}})_{{\mathbf m},{\mathbf k}\in\mathbb Z^n}$ and
$(g_{\mathbf{k l}})_{{\mathbf k},{\mathbf l}\in\mathbb Z^n}$ are
inverses of each other where
\begin{subequations}\label{fg}
\begin{gather}\label{fnk2}
f_{\mathbf{m k}}=\prod_{i=1}^n c_i(k_i)^{-1}
\prod_{1\le i<j\le n}\big(c_i(k_i)-c_j(k_j)\big)^{-1}\nonumber\\
\phantom{f_{\mathbf{m k}}=}{}\times\! \det_{1\le i,j\le
n}\!\!\Bigg[c_i(m_i)^{n+1-j} -a_i(m_i)^{n+1-j}
\frac{\big(c_i(m_i)-b/\prod\limits_{s=1}^n \! c_s(k_s)\big)}
{\big(a_i(m_i)-b/\prod\limits_{s=1}^n\! c_s(k_s)\big)}
\prod\limits_{s=1}^n\frac{\big(c_i(m_i)-c_s(k_s)\big)}
{\big(a_i(m_i)-c_s(k_s)\big)}\Bigg]\nonumber\\
\phantom{f_{\mathbf{m k}}=}{}\times
\prod_{i=1}^n\prod_{y_i=k_i+1}^{m_i}\Bigg[
\frac{\big(a_i(y_i)-b/\prod\limits_{j=1}^n c_j(k_j)\big)}
{\big(c_i(y_i)-b/\prod\limits_{j=1}^n c_j(k_j)\big)}
\prod_{j=1}^n\frac{\big(a_i(y_i)-c_j(k_j)\big)}
{\big(c_i(y_i)-c_j(k_j)\big)}\Bigg],
\end{gather}
and
\begin{gather}
g_{\mathbf{k l}}=\prod_{i=1}^n
\frac{\prod\limits_{y_i=l_i+1}^{k_i}\left[\big(a_i(y_i)-b/\prod\limits_{j=1}^n
c_j(k_j)\big)
\prod\limits_{j=1}^n\big(a_i(y_i)-c_j(k_j)\big)\right]}
{\prod\limits_{y_i=l_i}^{k_i-1}\left[\big(c_i(y_i)-b/\prod\limits_{j=1}^n
c_j(k_j)\big)
\prod\limits_{j=1}^n\big(c_i(y_i)-c_j(k_j)\big)\right]}.
\end{gather}
\end{subequations}
\end{Proposition}

\begin{Remark}\rm
The above Proposition generalizes Krattenthaler's matrix inverse
\cite{Kr} which is obtained when $n=1$. Note that for $n=1$ the
determinant in \eqref{fnk2} reduces (after relabeling) to
\begin{gather*}
c_m-a_m\frac{(c_m-b/c_k)}{(a_m-b/c_k)}\frac{(c_m-c_k)}{(a_m-c_k)}=
c_m\frac{(a_m-b/c_m)}{(a_m-b/c_k)}\frac{(a_m-c_m)}{(a_m-c_k)}
\end{gather*}
and the matrix inverse in \eqref{fg} (after relabeling) becomes
\begin{subequations}\label{fg1}
\begin{gather}\label{fnk1}
f_{m k}=\frac{(b-a_mc_m)(a_m-c_m)}{(b-a_kc_k)(a_k-c_k)}
\frac{\prod\limits_{y=k}^{m-1}(a_y-b/c_k)(a_y-c_k)}
{\prod\limits_{y=k+1}^{m}(c_y-b/c_k)(c_y-c_k)},
\\
\label{gkl1} g_{k
l}=\frac{\prod\limits_{y=l+1}^{k}(a_y-b/c_k)(a_y-c_k)}
{\prod\limits_{y=l}^{k-1}(c_y-b/c_k)(c_y-c_k)}.
\end{gather}
\end{subequations}
It is not dif\/f\/icult to see that this matrix inverse is
actually equivalent to its $b\to\infty$ special case:
\begin{subequations}\label{fg2}
\begin{gather}
f_{m k}=\frac{\prod\limits_{y=k}^{m-1}(a_y-c_k)}
{\prod\limits_{y=k+1}^{m}(c_y-c_k)},
\\
\label{gkl2} g_{k l}=\frac{(a_l-c_l)}{(a_k-c_k)}
\frac{\prod\limits_{y=l+1}^{k}(a_y-c_k)}
{\prod\limits_{y=l}^{k-1}(c_y-c_k)}.
\end{gather}
\end{subequations}
To recover \eqref{fg1} from \eqref{fg2}, do the substitutions
$a_y\mapsto a_y+b/a_y$, $c_y\mapsto c_y+b/c_y$, transfer some
simple factors from one matrix to the other, and simplify.

Other multidimensional generalizations of Krattenthaler's matrix
inverse were derived in \cite[Th.~3.1]{Sc1}, \cite{KS}, and
\cite{Sc2}.
\end{Remark}

In \cite[Th.~2.3]{LS} a matrix inverse slightly dif\/ferent to
Proposition~\ref{thmi2}, but equivalent to the latter, was given,
which contains a determinant in the entries of $(g_{\mathbf{k
l}})_{{\mathbf k},{\mathbf l}\in\mathbb Z^n}$ instead of
$(f_{\mathbf{m k}})_{{\mathbf m},{\mathbf k}\in\mathbb Z^n}$.

In the following we list a special case of Proposition~\ref{thmi2}
which is derived in a procedure analogous to the derivation of
\eqref{fg1} from \eqref{fg2}. This result is itself a
multidimensional generalization of the matrix inverse in
\eqref{fg1}.

\begin{Corollary}\label{thmi2c}
Let $b$ be an indeterminate and $a_i(k)$, $c_i(k)$ $(k\in \mathbb
Z,\  1\le i \le n)$ be arbitrary sequences of indeterminates. Then
the infinite lower-triangular $n$-dimensional matrices
$(f_{\mathbf{m k}})_{{\mathbf m},{\mathbf k}\in\mathbb Z^n}$ and
$(g_{\mathbf{k l}})_{{\mathbf k},{\mathbf l}\in\mathbb Z^n}$ are
inverses of each other where
\begin{subequations}
\begin{gather}
f_{\mathbf{m k}}=\prod_{i=1}^n \frac{c_i(m_i)^n}{c_i(k_i)^n}\,
\big(c_i(k_i)+b/c_i(k_i)\big)^{-1} \prod_{1\le i<j\le
n}\big[\big(1-b/c_i(k_i)c_j(k_j)\big)
\big(c_i(k_i)-c_j(k_j)\big)\big]^{-1}\nonumber\\
\phantom{f_{\mathbf{m k}}=}{}\times \det_{1\le i,j\le
n}\!\Bigg[\big(c_i(m_i)+b/c_i(m_i)\big)^{n+1-j}
-\big(a_i(m_i)+b/a_i(m_i)\big)^{n+1-j}\nonumber\\
\phantom{f_{\mathbf{m k}}=}{}\times
\prod_{s=1}^n\frac{\big(1-b/c_i(m_i)c_s(k_s)\big)\big(c_i(m_i)-c_s(k_s)\big)}
{\big(1-b/a_i(m_i)c_s(k_s)\big)\big(a_i(m_i)-c_s(k_s)\big)}\Bigg]\nonumber\\
\phantom{f_{\mathbf{m k}}=}{}\times
\prod_{i,j=1}^n\prod_{y_i=k_i+1}^{m_i}\Bigg[
\frac{\big(a_i(y_i)-b/c_j(k_j)\big)\big(a_i(y_i)-c_j(k_j)\big)}
{\big(c_i(y_i)-b/c_j(k_j)\big)\big(c_i(y_i)-c_j(k_j)\big)}\Bigg],
\end{gather}
and
\begin{gather}
g_{\mathbf{k l}}=\prod_{i,j=1}^n
\frac{\prod\limits_{y_i=l_i+1}^{k_i}\left[\big(a_i(y_i)-b/c_j(k_j)\big)
\big(a_i(y_i)-c_j(k_j)\big)\right]}
{\prod\limits_{y_i=l_i}^{k_i-1}\left[\big(c_i(y_i)-b/c_j(k_j)\big)
\big(c_i(y_i)-c_j(k_j)\big)\right]}.
\end{gather}
\end{subequations}
\end{Corollary}

\begin{proof}
In Proposition~\ref{thmi2}, f\/irst let $b\to\infty$, then perform
the substitutions $a_i(y_i)\mapsto a_i(y_i)+b/a_i(y_i)$ and
$c_i(y_i)\mapsto c_i(y_i)+b/c_i(y_i)$, for $1\le i\le n$. Finally,
transfer some factors from one matrix to the other.
\end{proof}

We give two important special cases of the above multidimensional
matrix inverses explicitly which involve the $q$-shifted
factorials $(a;q)_k$ def\/ined in \eqref{ipr}. For $n=1$ the
matrix inversions in Corollaries~\ref{cormi2} and \ref{cormi2c}
both reduce to Bressoud's~\cite{Br} matrix inverse, which he
directly extracted from the terminating very-well-poised
$_6\phi_5$ summation \eqref{sum65}.

\begin{Corollary}\label{cormi2}
Let $t_0,t_1,\dots,t_n$ and $u_1,\dots,u_n$ be indeterminates.
Then the infinite lower-triangular $n$-dimensional matrices
$(f_{\mathbf{m k}})_{{\mathbf m},{\mathbf k}\in\mathbb Z^n}$ and
$(g_{\mathbf{k l}})_{{\mathbf k},{\mathbf l}\in\mathbb Z^n}$ are
inverses of each other where
\begin{subequations}\label{cfg}
\begin{gather}
f_{\mathbf{m k}}=q^{(n-1)(|{\mathbf k}|-|{\mathbf m}|)}
\prod_{i=1}^nt_i^{n(m_i-k_i)}
\prod_{1\le i<j\le n}\big(u_iq^{k_i}-u_jq^{k_j}\big)^{-1}\nonumber\\
\phantom{f_{\mathbf{m k}}=}{}\times \det_{1\le i,j\le
n}\!\Bigg[\big(u_iq^{m_i}\big)^{n-j}\Bigg(1- t_i^{j-n-1}
\frac{\big(1-t_0u_iq^{m_i+|{\mathbf k}|}\big)}
{\big(1-t_0u_iq^{m_i+|{\mathbf k}|}/t_i\big)}
\prod_{s=1}^n\frac{\big(u_iq^{m_i}-u_sq^{k_s}\big)}
{\big(u_iq^{m_i}/t_i-u_sq^{k_s}\big)}\Bigg)\Bigg]\nonumber\\
\phantom{f_{\mathbf{m k}}=}{}\times \prod_{i=1}^n
\frac{(t_0u_iq^{1+k_i+|{\mathbf k}|}/t_i;q)_{m_i-k_i}}
{(t_0u_iq^{1+k_i+|{\mathbf k}|};q)_{m_i-k_i}} \prod_{i,j=1}^n
\frac{(q^{1+k_i-k_j}u_i/t_iu_j;q)_{m_i-k_i}}
{(q^{1+k_i-k_j}u_i/u_j;q)_{m_i-k_i}},\label{cormif}
\end{gather}
and
\begin{gather}
g_{\mathbf{k l}}= \prod_{i,j=1}^n
\frac{(qu_i/u_j;q)_{k_i-k_j}\,(t_ju_i/u_j;q)_{l_i-l_j}}
{(t_ju_i/u_j;q)_{k_i-k_j}\,(qu_i/u_j;q)_{l_i-l_j}}\nonumber\\
\phantom{g_{\mathbf{k l}}=}{}\times
\prod_{i=1}^n\frac{(t_0u_iq^{1+l_i+|{\mathbf
k}|}/t_i;q)_{k_i-l_i}} {(t_0u_iq^{l_i+|{\mathbf k}|};q)_{k_i-l_i}}
\prod_{i,j=1}^n \frac{(q^{l_i-l_j}t_ju_i/u_j;q)_{k_i-l_i}}
{(q^{1+l_i-l_j}u_i/u_j;q)_{k_i-l_i}}.\label{cormig}
\end{gather}
\end{subequations}
\end{Corollary}

\begin{Remark}\rm
The $t_i=t$, $i=0,1,\dots,n$ case of Corollary~\ref{cormi2} is
equivalent to Theorem 2.7 of~\cite{LS}.
\end{Remark}

\begin{proof}[Proof of Corollary~\ref{cormi2}]
We specialize Proposition~\ref{thmi2} by letting $b\mapsto
t_0^{-1}\prod\limits_{j=1}^nu_j$, $a_i(y_i)\mapsto
u_iq^{y_i}/t_i$, and $c_i(y_i)\mapsto u_iq^{y_i}$, for $1\le i\le
n$, and rewrite the expressions employing $q$-shifted factorial
notation. We f\/inally multiply $f_{\mathbf{m k}}$ by
$\prod\limits_{i=1}^n(q/t_i)^{n(k_i-m_i)}$ and $g_{\mathbf{k l}}$
by $\prod\limits_{i=1}^n(q/t_i)^{n(l_i-k_i)}$.
\end{proof}

\begin{Corollary}\label{cormi2c}
Let $t_0,t_1,\dots,t_n$ and $u_1,\dots,u_n$ be indeterminates.
Then the infinite lower-triangular $n$-dimensional matrices
$(f_{\mathbf{m k}})_{{\mathbf m},{\mathbf k}\in\mathbb Z^n}$ and
$(g_{\mathbf{k l}})_{{\mathbf k},{\mathbf l}\in\mathbb Z^n}$ are
inverses of each other where
\begin{subequations}\label{cfgc}
\begin{gather}
f_{\mathbf{m k}}=
\prod_{i=1}^nt_i^{n(m_i-k_i)}(u_iq^{k_i}+q^{-k_i}/au_i)^{-1}
\prod_{1\le i<j\le n}\big[\big(u_iq^{k_i}-u_jq^{k_j}\big)
\big(1-q^{-k_i-k_j}/au_iu_j\big)\big]^{-1}\nonumber\\
\phantom{f_{\mathbf{m k}}=}{}\times \det_{1\le i,j\le
n}\!\Bigg[\big(u_iq^{m_i}+q^{-m_i}/au_i\big)^{n+1-j}-
\big(u_iq^{m_i}/t_i+t_iq^{-m_i}/au_i\big)^{n+1-j}\nonumber\\
\phantom{f_{\mathbf{m k}}=}{}\times
\prod_{s=1}^n\frac{\big(1-q^{-m_i-k_s}/au_iu_s\big)
\big(u_iq^{m_i}-u_sq^{k_s}\big)}
{\big(1-t_iq^{-m_i-k_s}/au_iu_s\big)
\big(u_iq^{m_i}/t_i-u_sq^{k_s}\big)}\Bigg]\nonumber\\
\phantom{f_{\mathbf{m k}}=}{}\times \prod_{i,j=1}^n
\frac{(q^{1+k_i-k_j}u_i/t_iu_j,au_iu_jq^{1+k_i+k_j}/t_i;q)_{m_i-k_i}}
{(q^{1+k_i-k_j}u_i/u_j,au_iu_jq^{1+k_i+k_j};q)_{m_i-k_i}},\label{cormifc}
\end{gather}
and
\begin{gather}
g_{\mathbf{k l}}= \prod_{i,j=1}^n
\frac{(qu_i/u_j;q)_{k_i-k_j}\,(au_iu_jq/t_j;q)_{k_i+k_j}\,
(t_ju_i/u_j;q)_{l_i-l_j}\,(au_iu_j;q)_{l_i+l_j}}
{(t_ju_i/u_j;q)_{k_i-k_j}\,(au_iu_j;q)_{k_i+k_j}\,
(qu_i/u_j;q)_{l_i-l_j}\,(au_iu_jq/t_j;q)_{l_i+l_j}}\nonumber\\
\phantom{g_{\mathbf{k l}}=}{}\times \prod_{i,j=1}^n
\frac{(q^{l_i-l_j}t_ju_i/u_j,au_iu_jq^{l_i+l_j};q)_{k_i-l_i}}
{(q^{1+l_i-l_j}u_i/u_j,au_iu_jq^{l_i+l_j}/t_j;q)_{k_i-l_i}}.\label{cormigc}
\end{gather}
\end{subequations}
\end{Corollary}

\begin{proof}
We specialize Theorem~\ref{thmi2c} by letting $b\mapsto 1/a$,
$a_i(y_i)\mapsto u_iq^{y_i}/t_i$, and $c_i(y_i)\mapsto
u_iq^{y_i}$, for $1\le i\le n$, and rewrite the expressions
employing $q$-shifted factorial notation. We f\/inally multiply
$f_{\mathbf{m k}}$ by $\prod\limits_{i=1}^n(q/t_i)^{n(k_i-m_i)}$
and $g_{\mathbf{k l}}$ by
$\prod\limits_{i=1}^n(q/t_i)^{n(l_i-k_i)}$.
\end{proof}

\end{appendix}

 \subsection*{Acknowledgements}

I would like to thank Michel Lassalle for getting me involved into
Macdonald polynomials (especially concerning the issues related to
matrix inversion and explicit expressions) and his encouragement.
I would also like to express my sincere gratitude to the
organizers of the ``Workshop on Jack, Hall--Littlewood and
Macdonald Polynomials'' (ICMS, Edinburgh, September 23--26, 2003)
for inviting me to participate in that very stimulating workshop.
Among them, I am especially indebted to Vadim Kuznetsov whose
interest in explicit formulae for Macdonald polynomials served as
an inspiration for the present work.

The author was partly supported by FWF Austrian Science Fund
grants P17563-N13, and S9607 (the second is part of the Austrian
National Research Network ``Analytic Combinatorics and
Probabilistic Number Theory'').

\pdfbookmark[1]{References}{ref}
\LastPageEnding


\begin{thebibliography}{99}

\footnotesize\itemsep=0pt

\bibitem{An}
Andrews G.E., $q$-series: their development and application in
analysis, number theory, combinatorics, physics and computer
algebra, {\it CBMS Regional Conference Lectures Series}, Vol.~66,
Amer.\ Math.\ Soc., Providence, RI, 1986.

\bibitem{AAR}
Andrews G.E., Askey  R., Roy R., Special functions, {\it
Encyclopedia of Mathematics and Its Applications}, Vol.~71,
Cambridge University Press, Cambridge, 1999.

\bibitem{BF}
 Baker T.H., Forrester P.J.,
Transformation formulas for multivariable basic hypergeometric
series, {\em Meth.\ Appl.\ Anal.} \textbf{6} (1999), 147--164,
\href{http://arxiv.org/abs/math.QA/9803146}{math.QA/9803146}.

\bibitem{Bh}
 Bhatnagar G.,
A multivariable view of $q$-series, in Special Functions and
Dif\/ferential Equations, Editors K.~Srinivasa Rao,
R.~Jagannathan, G.~Vanden Berghe and J.\ Van der Jeugt,
Proceedings of a Workshop, WSSF'97 (January 13--24, 1997, Madras,
India),
 Allied Publ., New Delhi, 1998, 25--30.

\bibitem{BM}
 Bhatnagar G., Milne S.C.,
Generalized bibasic hypergeometric series and their $U(n)$
extensions, {\em Adv. Math.} \textbf{131} (1997), 188--252.

\bibitem{Br}
 Bressoud D.M.,
A matrix inverse, {\em Proc.\ Amer.\ Math.\ Soc.} \textbf{88}
(1983), 446--448.

\bibitem{Bw}  Bromwich T.J.l'A.,
An introduction to the theory of inf\/inite series, 2nd ed.,
Macmillan, London, 1949.

\bibitem{C}
 Cauchy A.-L.,
M\'emoire sur les fonctions dont plusieurs valeurs sont li\'ees
entre elles par une \'equation lin\'eaire, et sur diverses
transformations de produits compos\'es d'un nombre ind\'ef\/ini de
facteurs, {\em C.\ R.\ Acad.\ Sci.\ Paris} {\bf XVII} (1843), 523;
{\em Oeuvres de Cauchy}, $1^{\text{re}}$ s\'erie, T.~VIII,
Gauthier-Villars, Paris, 1893, 42--50.

\bibitem{DK}
 Danilov V., Koshevoy G.,
Continuous combinatorics, Preprint, 2005, 12 pages.

\bibitem{DG}
 Denis R.Y., Gustafson R.A.,
An $SU(n)$ $q$-beta integral transformation and multiple
hypergeometric series identities, {\em SIAM J.\ Math.\ Anal.}
\textbf{23} (1992), 552--561.

\bibitem{FT}
 Frenkel I.B., Turaev V.G., Elliptic solutions of the
Yang--Baxter equation and modular hypergeometric functions, in The
Arnold--Gelfand Mathematical Seminars, Editors V.I.~Arnold,
I.M.~Gelfand, V.S.~Retakh and M.~Smirnov, Birkh\"auser, Boston,
1997, 171--204.

\bibitem{GR}
 Gasper G., Rahman M.,
Basic hypergeometric series,  2nd ed., {\it Encyclopedia of
Mathematics and Its Applications}, Vol.~96, Cambridge University
Press, Cambridge, 2004.

\bibitem{G1}
Gustafson R.A., Multilateral summation theorems for ordinary and
basic hypergeometric series in $U(n)$, {\em SIAM J.\ Math.\ Anal.}
\textbf{18} (1987), 1576--1596.

\bibitem{G2}
 Gustafson R.A.,
The Macdonald identities for af\/f\/ine root systems of classical
type and hypergeometric series very well-poised on semi-simple Lie
algebras, in Ramanujan International Symposium on Analysis
(December 26--28, 1987, Pune, India), Editor N.K.~Thakare, 1989,
187--224.

\bibitem{H}
 Heine E., Untersuchungen \"uber die Reihe $\ldots$,
{\em J.\ Reine Angew.\ Math.} \textbf{34} (1847), 285--328.

\bibitem{Ho}
Holman W.J. III, Summation theorems for hypergeometric series in
$U(n)$, {\em SIAM J.\ Math.\ Anal.} \textbf{11} (1980), 523--532.

\bibitem{HBL}
 Holman W.J. III,  Biedenharn L.C., Louck J.D.,
On hypergeometric series well-poised in $SU(n)$, {\em SIAM J.\
Math.\ Anal.} \textbf{7} (1976), 529--541.

\bibitem{Is}
 Ismail M.E.H.,
Classical and quantum orthogonal polynomials in one variable, {\it
Encyclopedia of Mathematics and its Applications}, Vol.~98,
Cambridge University Press, Cambridge, 2005.

\bibitem{Ja}
 Jackson F.H., Summation of $q$-hypergeometric series,
{\em Messenger of Math.} {\bf 57} (1921), 101--112.

\bibitem{JJ}
 Jing N.H., J\'{o}zef\/iak T.,
A formula for two-row Macdonald functions, {\em Duke Math.\ J.}
\textbf{67} (1992), 377--385.

\bibitem{Ka}
 Kadell K.,
The Schur functions for partitions with complex parts, {\em
Contemp.\ Math.} \textbf{254} (2000), 247--270.

\bibitem{KN}
 Kajihara Y., Noumi M.,
Raising operators of row type for Macdonald polynomials, {\em
Compos.\ Math.} \textbf{120} (2000), 119--136,
\href{http://arxiv.org/abs/math.QA/9803151}{math.QA/9803151}.

\bibitem{Kan1}
 Kaneko J.,
$q$-Selberg integrals and Macdonald polynomials, {\em Ann.\ Sci.\
\'Ecole Norm.\ Sup.} \textbf{29} (1996), 583--637.

\bibitem{Kan2}
 Kaneko J.,
A ${}_1\Psi_1$ summation theorem for Macdonald polynomials, {\em
Ramanujan J.\ }\textbf{2} (1998), 379--386.

\bibitem{Ko}
 Koornwinder T.H.,
Self-duality for $q$-ultraspherical polynomials associated with
root system $A_n$, unpublished handwritten manuscript,
1988, 17 pages.\\
Available at
\url{http://remote.science.uva.nl/~thk/art/informal/dualmacdonald.pdf}.

\bibitem{Kr}
 Krattenthaler C., A new matrix inverse,
{\em Proc.\ Amer.\ Math.\ Soc.} \textbf{124} (1996), 47--59.

\bibitem{KS}
 Krattenthaler C., Schlosser M.J.,
A new multidimensional matrix inverse with applications to
multiple $q$-series, {\em Discrete Math.} \textbf {204} (1999),
249--279.

\bibitem{La}
 Lassalle M.,
Explicitation des polyn\^omes de Jack et de Macdonald en longueur
trois, {\em C.\ R.\ Acad.\ Sci.\ Paris S\'er.\ I Math.}
\textbf{333} (2001), 505--508.

\bibitem{La1}
 Lassalle M.,
Une $q$-sp\'ecialisation pour les fonctions sym\'etriques
monomiales, {\em Adv.\ Math.} \textbf{162} (2001), 217--242,
\href{http://arxiv.org/abs/math.CO/0004019}{math.CO/0004019}.

\bibitem{La2}
 Lassalle M.,
A short proof of generalized Jacobi--Trudi expansions for
Macdonald polynomials, {\em Contemp.\ Math.} \textbf{417} (2006),
271--280,
\href{http://arxiv.org/abs/math.CO/0401032}{math.CO/0401032}.

\bibitem{LS}
 Lassalle M., Schlosser M.J.,
Inversion of the Pieri formula for Macdonald polynomials, {\em
Adv. Math.} \textbf{202} (2006), 289--325,
\href{http://arxiv.org/abs/math.CO/0402127}{math.CO/0402127}.

\bibitem{Ma1}
 Macdonald I.G., A new class of symmetric functions,
\emph{S\'em.\ Lothar.\ Combin.} \textbf{20} (1988), Art.~B20a,
41~pages.

\bibitem{Ma}
 Macdonald I.G., Symmetric functions and hall polynomials, 2nd ed.,
Clarendon Press, Oxford, 1995.

\bibitem{Ma2}
 Macdonald I.G., Symmetric functions and orthogonal polynomials,
Dean Jacqueline B.~Lewis Memorial Lectures presented at Rutgers
University, New Brunswick, NJ, \emph{University Lecture Series},
Vol.~12, Amer.\ Math.\ Soc., Providence, RI, 1998.

\bibitem{Mi1}
 Milne S.C.,
An elementary proof of the Macdonald identities for
$A_\ell^{(1)}$, {\em Adv.\ Math.} \textbf{57} (1985), 34--70.

\bibitem{Mi2}
 Milne S.C.,
Basic hypergeometric series very well-poised in $U(n)$, {\em J.\
Math.\ Anal.\ Appl.} \textbf{122} (1987), 223--256.

\bibitem{Mi7}
 Milne S.C.,
Multiple $q$-series and $U(n)$ generalizations of Ramanujan's
$_1\psi_1$ sum, in Ramanujan Revisited, Editors G.E.~Andrews et
al., Academic Press, New York, 1988, 473--524.

\bibitem{Mi3}
 Milne S.C.,
The multidimensional $_1\Psi_1$ sum and Macdonald identities for
$A_l^{(1)}$, in Theta Functions Bowdoin (1987), Editors
L.~Ehrenpreis and R.C.~Gunning, {\em Proc.\ Sympos.\ Pure Math.}
\textbf{49} (1989), 323--359.

\bibitem{Mi4}
 Milne S.C.,
A $q$-analog of a Whipple's transformation for hypergeometric
series in $U(n)$, {\em Adv.\ Math.} \textbf{108} (1994), 1--76.

\bibitem{Mi5}
 Milne S.C.,
Balanced $_3\phi_2$ summation theorems for $U(n)$ basic
hypergeometric series, {\em Adv.\ Math.} \textbf{131} (1997),
93--187.

\bibitem{Mi6}
 Milne S.C.,
Transformations of $U(n+1)$ multiple basic hypergeometric series,
in  Physics and Combinatorics, Editors A.N.~Kirillov, A.~Tsuchiya
and H.~Umemura, Proceedings of the Nagoya 1999 International
Workshop (August 23--27, 1999, Nagoya University, Japan), World
Scientif\/ic, Singapore, 2001, 201--243.

\bibitem{ML}
 Milne S.C., Lilly G.M.,
Consequences of the $A_l$ and $C_l$ Bailey transform and Bailey
lemma, {\em Discrete Math.} \textbf{139} (1995), 319--346,
\href{http://arxiv.org/abs/math.CA/9204236}{math.CA/9204236}.

\bibitem{MN}
 Milne S.C.,  Newcomb J.W.,
$U(n)$ very-well-poised $_{10}\phi_9$ transformations, {\em J.\
Comput.\ Appl.\ Math.} \textbf{68} (1996), 239--285.

\bibitem{Ok}
 Okounkov A., $BC$-type interpolation Macdonald polynomials and
binomial formula for Koornwinder polynomials, \emph{Transform.\
Groups} \textbf{3} (1998), 181--207,
\href{http://arxiv.org/abs/q-alg/9611011}{q-alg/9611011}.

\bibitem{Ra1}
 Rains E.M., $BC_n$-symmetric polynomials,
\emph{Transform.\ Groups} \textbf{10} (2005), 63--132,
\href{http://arxiv.org/abs/math.QA/0112035}{math.QA/0112035}.

\bibitem{Ra2}
 Rains E.M., $BC_n$-symmetric Abelian functions,
\emph{Duke Math.\ J.} \textbf{135} (2006), 99--180,
\href{http://arxiv.org/abs/math.CO/0402113}{math.CO/0402113}.

\bibitem{Ro}
 Rogers R.J.,
Third memoir on the expansion of certain inf\/inite products, {\em
Proc.\ London Math.\ Soc.} \textbf{26} (1894), 15--32.

\bibitem{R1}
 Rosengren H.,
Elliptic hypergeometric series on root systems, {\em Adv.\ Math.}
\textbf{181} (2004), 417--447,
\href{http://arxiv.org/abs/math.CA/0207046}{math.CA/0207046}.

\bibitem{R2}
 Rosengren H.,
Reduction formulas for Karlsson--Minton-type hypergeometric
functions, {\em Constr.\ Approx.} \textbf{20} (2004), 525--548,
\href{http://arxiv.org/abs/math.CA/0202232}{math.CA/0202232}.

\bibitem{Sc1}
 Schlosser M.J.,
Multidimensional matrix inversions and $A_r$ and $D_r$ basic
hypergeometric series, {\em Rama\-nujan~J.} \textbf{1} (1997),
243--274.

\bibitem{Sc2}
 Schlosser M.J.,
A new multidimensional matrix inversion in $A_r$, {\em Contemp.\
Math.} \textbf{254} (2000), 413--432.

\bibitem{Sp}
 Spiridonov V.P.,
Elliptic hypergeometric functions, Dr.Sc. Thesis,
JINR, Dubna, Russia, 2004, 218 pages.

\bibitem{St}
 Stanton D.,
An elementary approach to the Macdonald identities, in $q$-Series
and Partitions, Editor D.~Stanton, {\it The IMA Volumes in
Mathematics and Its Applications}, Vol.~18, Springer-Verlag, 1989,
139--150.

\bibitem{W}
 Warnaar S.O.,
$q$-Selberg integrals and Macdonald polynomials, {\em Ramanujan
J.} \textbf{10} (2005), 237--268.

\bibitem{Wt}
 Watson G.N.,
A treatise on the theory of Bessel functions, 2nd ed., Cambridge
University Press, Cambridge, 1966.

\end{thebibliography}
\end{document}